\documentclass[12pt]{article}
\usepackage{amsmath,amssymb,amsfonts,amsthm}
\usepackage[dvips]{graphicx}

\def\UseSection{
        \numberwithin{equation}{section}
    \theoremstyle{plain}
        \newtheorem{theorem}    {Theorem}[section]
        \DefineTheorems 
}

\def\DefineTheorems{
    
    \newtheorem{lemma}      [theorem] {Lemma}
    
    \newtheorem{prop}       [theorem] {Proposition}
    
    \newtheorem{cor}        [theorem] {Corollary}

    \theoremstyle{definition}
    \newtheorem{defn}       [theorem] {Definition}

    \theoremstyle{definition}

}

\newcommand{\bt}   {\begin{theorem}}
\newcommand{\et}   {\end  {theorem}}
\newcommand{\bl}   {\begin{lemma}}
\newcommand{\el}   {\end  {lemma}}
\newcommand{\bp}   {\begin{prop}}
\newcommand{\ep}   {\end  {prop}}
\newcommand{\bc}   {\begin{cor}}
\newcommand{\ec}   {\end  {cor}}
\newcommand{\bd}   {\begin{defn}}
\newcommand{\ed}   {\end  {defn}}

\newcommand{\ba}   {\begin{array}}
\newcommand{\ea}   {\end  {array}}
\newcommand{\be}   {\begin{enumerate}}
\newcommand{\ee}   {\end  {enumerate}}
\newcommand{\bi}   {\begin{itemize}}
\newcommand{\ei}   {\end  {itemize}}

\def\eq#1\en{\begin{equation}#1\end{equation}}
\def\eqsplit#1\ensplit{
    \begin{equation}\begin{split}#1\end{split}\end{equation}
    }
\def\eqalign#1\enalign{
    \begin{align}#1\end{align}
    }
\def\eqmul#1\enmul{
    \begin{multline}#1\end{multline}
    }
\newcommand{\eqarrstar} {\begin{eqnarray*}}
\newcommand{\enarrstar} {\end{eqnarray*}}
\newcommand{\eqarray}   {\begin{eqnarray}}
\newcommand{\enarray}   {\end{eqnarray}}
\newcommand{\nnb}   {\nonumber \\}

\newcommand{\lbeq}[1]  {\label{e:#1}}
\newcommand{\refeq}[1] {\eqref{e:#1}}    

\makeatletter
\newcommand{\labelcounter}[2]{{%
    \stepcounter{#1}
    \protected@write\@auxout{}%
    {\string\newlabel{#2}{{\csname the#1\endcsname}{\thepage}}}%
    {\ref{#2}}
    }}
\makeatother
\newcommand{\Nbold} {{\mathbb N}}

\newcommand{\Pbold} {{\mathbb P}}

\newcommand{\Rbold} {{\mathbb R}}

\newcommand{\Zbold} {{\mathbb Z}}

\newcommand{\Ccal}   {\mathcal{C}}

\newcommand{\Ocal}   {\mathcal{O}}
\newcommand{\Pcal}   {\mathcal{P}}

\newcommand{\Rcal}   {\mathcal{R}}

\newcommand{\Rd}    {{ {\Rbold}^d}}
\newcommand{\Zd}    {{ {\Zbold}^d }}

\newcommand{\spose}[1] {{\hbox to 0pt{#1\hss}} }
\newcommand{\ltapprox} {\mathrel{\spose{\lower 3pt\hbox{$\mathchar"218$}}
 \raise 2.0pt\hbox{$\mathchar"13C$}}}
\newcommand{\gtapprox} {\mathrel{\spose{\lower 3pt\hbox{$\mathchar"218$}}
 \raise 2.0pt\hbox{$\mathchar"13E$}}}

\UseSection  
\renewcommand{\to}      {\rightarrow}

\newcounter{countC}  
\newcounter{countR}  
\newcommand{\sumtwo}[2]{\sum_{ \mbox{ \scriptsize
    $\begin{array}{c}
                        {#1} \\ {#2}
                        \end{array} $ }
    }
}

\newcommand{\R}{\Rbold}
\newcommand{\Z}{\Zbold}
\newcommand{\N}{\Nbold}
\newcommand{\conn}{\to}

\newcommand{\dbc}{\Rightarrow}

\newcommand{\nn}{\nonumber}

\title  {
        A generalised inductive approach 
        to the lace expansion
        }

\author{\begin{tabular}{ccc}
    Remco van der Hofstad           & Gordon Slade \\
        Stieltjes Institute for Mathematics & Department of Mathematics \\
        Delft University            & University of British Columbia \\
        Mekelweg 4              & Vancouver, BC V6T 1Z2 \\
        2628 CD Delft, The Netherlands      & Canada \\
        {\tt R.W.vanderHofstad@its.tudelft.nl}  & {\tt slade@math.ubc.ca}\\
        \end{tabular}}

\date{August 31, 2000}

\begin{document}

\maketitle


\begin{abstract}
The lace expansion is a powerful tool for analysing the critical
behaviour of self-avoiding walks and percolation.  It gives rise
to a recursion relation which we abstract and study using an
adaptation of the inductive method introduced by den Hollander and
the authors. We give conditions under which the solution to the
recursion relation behaves as a Gaussian, both in Fourier space
and in terms of a local central limit theorem.
These conditions are
shown elsewhere to hold for sufficiently spread-out models of
networks of self-avoiding walks in dimensions $d>4$, and for critical oriented
percolation in
dimensions $d+1>5$, providing a unified approach and an essential
ingredient for
a detailed analysis of the branching behaviour of
these models.
\end{abstract}


\section{The lace expansion}
\label{sec-lace}

\subsection{The recursion relation}
\label{sec-recrel}

When applied to self-avoiding walks or oriented percolation, the
lace expansion gives rise to a convolution recursion relation of
the form
        \eq
        \lbeq{fkrec}
        f_{n+1}(k;z) =
        \sum_{m=1}^{n+1} g_m(k;z)f_{n+1-m}(k;z) + e_{n+1}(k;z)
        \quad \quad (n \geq 0),
        \en
with $f_0(k;z) = 1$.
Here, $k \in [-\pi,\pi]^d$ is a parameter dual to a
spatial lattice variable $x \in \Zd$.  The positive parameter $z$ is
related to the activity for self-avoiding walks, and to
the bond occupation probability for oriented percolation.

Such a recursion arises, in a simple form, already for ordinary random walk.
In fact, let $D$ be a non-negative function on $\Zd$ which respects
the lattice symmetries
and obeys $\sum_{x \in \Zd}D(x)=1$.  Let $z>0$, and let
$q_n(x;z)$ denote the $n$-step transition function of
the random walk on $\Zd$ with $1$-step transition function
$D(x)$, multiplied by $z^n$.   Then $q_1(x;z)=zD(x)$.
Set $q_0(x;z) = \delta_{0,x}$.
Then $q_n$ obeys
\eq
\lbeq{qnrec}
    q_n(x;z) = \sum_{y \in \Zd} zD(y) q_{n-1}(x-y;z) \quad (n \geq 1).
\en
Given an absolutely summable function $h$ on $\Zd$, we denote its Fourier
transform by
\eq
    \hat{h}(k) = \sum_{x \in \Zd} h(x) e^{ik\cdot x}
    \quad (k \in [-\pi,\pi]^d).
\en
Then
\refeq{qnrec} implies that
\eq
\lbeq{qnkrec}
    \hat{q}_n(k;z) = z\hat{D}(k) \hat{q}_{n-1}(k;z).
\en
The recursion relation \refeq{qnkrec} is in the form \refeq{fkrec}, with
$f_n(k;z) = \hat{q}_n(k;z)$, $g_1(k;z) = z\hat{D}(k)$, $g_m = 0$
for $m \geq 2$, and $e_n = 0 $ for $n \geq 1$.
It is readily solved to give $\hat{q}_n(k;z) = z^n \hat{D}(k)^n$.
Assuming that $D$ has a finite second moment $\sigma^2 = \sum_x |x|^2D(x)$,
Taylor expansion of $\hat{D}$ then gives the central limit theorem
\eq
\lbeq{qnclt}
    \lim_{n \to \infty}
    \hat{q}_n \big( \frac{k}{\sqrt{\sigma^2 n}} ; 1 \big)
    = e^{-k^2/2d}.
\en
The value $z=1$ plays the role of a critical point, with solutions
of the recursion exponentially growing for $z>1$ and exponentially decaying
if $z<1$.

Our goal here is to prove that the solutions of
certain more elaborate recursion
relations still exhibit Gaussian behaviour, as in \refeq{qnclt}.
The recursion relations
of interest arise, as we will explain below in Section~\ref{sec-examples},
in applying the lace expansion to
self-avoiding walks or oriented percolation.
The lace
expansion produces a recursion relation of the form \refeq{fkrec}, with
$g_1(k;z) = z\hat{D}(k)$, and with
explicit expressions for the functions $g_m$ and $e_n$ for $m \geq 2$
and $n \geq 1$.
For oriented percolation, $g_n$ and $e_n$ are almost identical, whereas
for self-avoiding walk $e_n=0$ for all $n$.
A small parameter is introduced into the models, with the result that
$g_m$ and $e_n$ are small for $m \geq 2$ and $n \geq 1$.  They are small
uniformly in $m$, as well as
in terms of their rapid decay as $m$ becomes large.
The recursion relation \refeq{fkrec} can thus be considered to be a small
perturbation of the simple equation \refeq{qnkrec}.
The recursion turns out to have
a critical value $z_c$, and, as we have seen above for
ordinary random walk, we will set $z=z_c$ to observe Gaussian
behaviour.

We analyse the recursion relation \refeq{fkrec} using
induction on $n$.
Our method is a generalisation of the inductive method
introduced by den Hollander and the authors \cite{HHS98} to study
weakly self-avoiding walk for $d>4$.  Here, we are considering long-range
or ``spread-out''
models.  This necessitates a number of modifications to the analysis
of \cite{HHS98},
as the small parameter is less explicit for the spread-out
model. In addition, the analysis of \cite{HHS98} has been
simplified in some respects.  Our choice of induction hypotheses
has been strongly motivated by those used in \cite{HHS98}, which
in turn were motivated by earlier work.
The inductive approach leads to
more detailed results than have been
obtained previously for self-avoiding walk or oriented
percolation using generating functions.
In a different vein, it has also been adapted to prove
ballistic behaviour for weakly self-avoiding walk for $d=1$ \cite{Hofs99}.

We will prove Gaussian behaviour of $f_n$ under certain
assumptions on $e_n$ and $g_n$, as well as some modest assumptions
on $f_n$.  Since the assumptions on $e_n$ and $g_n$ are motivated by
the lace expansion, we will give a brief review of relevant aspects
of the lace expansion, in Section~\ref{sec-examples}, to provide some
necessary context.
Our assumptions are shown elsewhere
to apply for sufficiently spread-out models of self-avoiding walks
on $\Zd$ \cite{HS01b} and critical oriented bond percolation on
$\Zd \times \Z_+$ \cite{HHS01b}, in dimensions $d>4$.
Our results provide an essential ingredient in both cases, and form
the basis for a proof of Gaussian behaviour of sufficiently spread-out
polymer networks for $d>4$, and for a proof that the moment measures
of sufficiently spread-out
critical oriented percolation converge to those of super-Brownian
motion for $d>4$.  In this paper, we extract an
important model-independent aspect of
the proof of these results.  The precise form of the functions $e_n$ and
$g_n$ is model-dependent, and the verification of the assumptions on
them does depend on the model.

\subsection{Assumptions on the recursion relation}
\label{sec-assumpt}

Before stating the assumptions, we first introduce some notation.
For $x=(x_1,\ldots,x_d)$ in $\Zd$ or $\Rd$, we write $\|x\|_p = (\sum_{i=1}^d
|x_i|^p)^{1/p}$ for $p \in [1,\infty)$, and $\|x\|_\infty =
\max_{i \in \{1,\ldots,d\}}|x_i|$.  For a function $h$ on $[-\pi,\pi]^d$,
we write $\|h\|_p = [(2\pi)^{-d}\int_{[-\pi,\pi]^d} |h(k)|^p d^dk ]^{1/p}$
and use $\|h\|_\infty$ to denote the essential supremum of $h$.
Similarly, for a function $H$ on $\Zd$, we write $\|H\|_p = [\sum_{x \in \Zd}
|H(x)|^p]^{1/p}$ and $\|H\|_\infty = \sup_x |H(x)|$.

Throughout this paper, we make the following four assumptions,
Assumptions~S, D, E, G, on the
quantities appearing in the recursion equation \refeq{fkrec}.

The first assumption,
Assumption S, requires that the functions appearing in \refeq{fkrec}
respect the lattice symmetries, and that $f_n$ remains bounded
in a weak sense.

\smallskip \noindent
{\bf Assumption S}. For every $n\in\N$ and $z>0$, the mapping
$k\mapsto f_n(k;z)$ is symmetric under replacement of any component
$k_i$ of $k$ by
$-k_i$, and under permutations of the components of $k$.  The same
holds for $e_n(\cdot;z)$ and $g_n(\cdot; z)$.  In addition, for
each $n$, $|f_n(k;z)|$ is bounded uniformly in $k \in
[-\pi,\pi]^d$ and $z$ in a neighbourhood of $1$ (which may depend on $n$).

\smallskip
Our next assumption incorporates a ``spread-out'' aspect to the recursion
equation.  It introduces a function
$D$ which defines the
underlying random walk model, about which \refeq{fkrec} is a perturbation.
The assumption involves a non-negative parameter $L$, which will
be taken to be large, and which serves to spread out the steps of the
random walk over a large set.
We write $D=D_L$ in the statement of the assumption to emphasise this
dependence, but the subscript will not be retained elsewhere.
An example of a family of $D$'s obeying
the assumption is given following its statement.  The assumption implies,
in particular, that $D$ has a finite second moment, and we define
\eq
\lbeq{sigdef}
    \sigma^2 = - \nabla^2 \hat{D}(0).
\en
Our assumptions will involve a parameter $d$, which corresponds to the
spatial dimension in our applications.  Our theorems will assume that
$d>4$.  Error terms will involve an arbitrarily
chosen $\delta \in (0, 1 \wedge \frac{d-4}{2})$, which should be thought
of as being close to $1 \wedge \frac{d-4}{2}$.
We will write
    \eq \lbeq{adef}
    a(k) = 1 - \hat{D}(k).
    \en

\smallskip \noindent
{\bf Assumption D}.
Let $D_L$ be a non-negative function on $\Zd$.
We assume that
    \eq
    f_1(k;z) = z \hat{D}_L(k), \quad e_1(k;z)=0.
    \en
In particular, this implies that
$g_1(k;z)=z\hat{D}_L(k)$.
The function $D_L$ is assumed to be normalised, so that
$\hat{D}_L(0) =1$, and to
have $2+2\epsilon$ moments, for some
$\epsilon > \delta$, i.e.,
\eq
    \lbeq{momentD}
    \sum_{x\in \Z^d} |x|^{2+2\epsilon} D_L(x) <\infty .
\en
We assume that there is a constant $C$ such that
\eq
\lbeq{beta, sigmadef}
    \|D_L\|_\infty \leq CL^{-d},\qquad \sigma^2 = \sigma^2_L
    \leq CL^2,
\en
for all $L \geq 1$.  In addition, we assume that there exist
constants $\eta,c_1,c_2 >0$ such that
\eq
\lbeq{Dbound1}
    c_1 L^2 k^2 \leq a_L(k) \leq c_2 L^2 k^2
    \quad
    (\|k\|_\infty \leq L^{-1}),
\en
\eq
\lbeq{Dbound2}
    a_L(k) > \eta  \quad (\|k\|_\infty \geq L^{-1}),
\en
\eq
\lbeq{Dbound3}
    a_L(k) < 2-\eta \quad (k \in [-\pi,\pi]^d).
\en

\smallskip \noindent
{\bf Example.}
Let $h$ be a non-negative bounded
function on $\Rd$ which is almost everywhere continuous and
symmetric under the lattice symmetries.  Assume that there is a function
$H$ on $\Rd$ with $H(te)$
non-increasing in $t \geq 0$ for every unit vector $e \in \Rd$,
such that $h(x) \leq H(x)$ for all $x \in \Rd$.
Assume that $\int_{\Rd} h(x) d^d x =1$ and
$\int_{\R^d} |x|^{2+2\epsilon} h(x) d^dx <\infty$ for some $\epsilon > \delta$.
The monotonicity hypothesis implies that
$\sum_{x \in \Zd} h(x/L) < \infty$ for all $L$.
We will verify in
Appendix~\ref{sec-app} that the function
    \eq
    \lbeq{Dbasic2}
    D_L(x) = \frac{h(x/L)}{\sum_{x\in \Z^d}h(x/L)}
    \en
obeys the conditions of Assumption~D, when $L$ is large enough.
A good example is given by the function $h(x)=2^{-d}$
for $x \in [-1,1]^d$, $h(x)=0$ otherwise.
In this case, $D_L$ is uniform on the cube $[-L,L]^d \cap \Zd$.


\medskip
We will take $L$ large, yielding a small parameter
\eq
    \beta = L^{-d}.
\en
Each of $e_n$, $f_n$ and $g_n$ depend on
$L$, but we do not make this dependence explicit in the notation.

The particular coefficients $e_m$ and $g_m$ arising in lace
expansion analyses have Feynman diagram representations, as we
discuss further in Section~\ref{sec-examples}.  These
Feynman diagrams can be bounded in terms of $f_j$, for $j<m$.
This {\em strict}\/ inequality  $j<m$
means that the right side of \refeq{fkrec}
can be studied in terms of $f_m$ for $m < n+1$.  This mechanism drives the
induction.

Our next two assumptions summarise the
manner in which bounds on $f_m$ imply bounds on $e_n$
and $g_n$.  Verification of these assumptions requires a
model-dependent analysis, which is carried out for
self-avoiding walks in \cite{HS01b} and for oriented percolation in
\cite{HHS01b}.  A brief indication of how this is done is given in
Section~\ref{sec-examples}.
However, given the assumptions, the remaining
inductive analysis, which we present in this paper, is
model-independent.

The relevant bounds on $f_m$, which {\em a priori}\/ may or may
not be satisfied, are
    \eq \lbeq{fbds}
    \|\hat{D}^2 f_m(\cdot;z)\|_1\leq K\beta m^{-d/2},
    \quad  | f_m(0;z)|\leq K, \quad
    |\nabla^2 f_m(0;z)|\leq K \sigma^2 m,
    \en
for some positive constant $K$.
The factor $\hat{D}(k)^2$ appearing in the first bound can be
understood from the fact that $\|\hat{D}^2\|_1$ is the probability
that a random walk taking steps with step distribution $D$ returns to
the origin after two steps, which is bounded by $\|D\|_\infty \leq
\Ocal(\beta)$.  Because of
this, the factor $\hat{D}(k)^2$ is helpful in extracting the
small factor $\beta$ appearing in the upper bound.

\smallskip \noindent
{\bf Assumption E}. Suppose that there is a $K>1$ and an interval
$I \subset [1-\alpha,1+\alpha]$, for some $\alpha \in (0,1)$, such that
the bounds \refeq{fbds} hold for $z \in I$ and for $m$ with $1
\leq m \leq n$.  Then, there exist $C_e=C_e(K)<\infty$ and $L_0>0$ such
that for all $L \geq L_0$, $z \in I$, $k \in [-\pi,\pi]^d$
and $2 \leq m\leq n+1$, the following bounds hold:
    $$
    |e_m(k;z)|\leq C_e(K) \beta m^{-d/2},
    \quad |e_m(k;z)-e_m(0;z)|\leq
    C_e(K) a(k) \beta m^{-(d-2)/2}.
    $$

\smallskip \noindent
{\bf Assumption G}. Suppose that there is a $K>1$ and an interval
$I \subset [1-\alpha,1+\alpha]$, for some $\alpha \in (0,1)$, such that
the bounds \refeq{fbds} hold for $z \in I$ and for $m$ with $1
\leq m \leq n$.  Then, there exist $C_g=C_g(K)<\infty$ and $L_0>0$ such
that for all $L \geq L_0$, $z \in I$, $k \in [-\pi,\pi]^d$
and $2 \leq m\leq n+1$, the following bounds hold:
    $$
    |g_m(k;z)|\leq C_g(K)\beta m^{-d/2},
    \quad |\nabla^2 g_m(0;z)|\leq
    C_g(K) \sigma^2 \beta m^{-(d-2)/2},
    $$
    $$ |\partial_z g_m(0;z)|\leq
    C_g(K)\beta m^{-(d-2)/2},
    $$
    $$|g_m(k;z)-g_m(0;z)- a(k) \sigma^{-2}
    \nabla^2 g_m(0;z)| \leq C_g(K)\beta
    a(k)^{1+\epsilon'}m^{-(d-2-2\epsilon')/2},
    $$
with the last bound
valid for any $\epsilon' \in [0,\epsilon]$, with $\epsilon$ given in
\refeq{momentD}.

\smallskip
We emphasize that it is the assumption that bounds on $f_m$ for $m
\leq n$ imply bounds on $g_m$ and $e_m$ for $m \leq n+1$ that
permits an inductive analysis.

\subsection{Main results}
\label{sec-thm}

We assume henceforth that Assumptions S, D, E and G all hold. Our
first theorem shows that there is a value $z_c$ for which $f_n$
behaves as a Gaussian.  Recall that $\sigma$ was defined in \refeq{sigdef}.

\begin{theorem}
\label{thm-1}
Let $d>4$ and fix an arbitrary $\delta \in (0,1\wedge \frac{d-4}{4})$.
There exist positive $L_0 = L_0(d)$, $z_c=z_c(d,L)$,
$A=A(d,L)$, and
$v = v(d,L)$,
such that for $L \geq L_0$,\\
(a)
        \eq
        f_n\Big(\frac{k}{\sqrt{v \sigma^2 n}};z_c\Big) =
        A e^{-\frac{k^2}{2d}}
        [1 + {\cal O}(k^2n^{-\delta})+{\cal O}(n^{-(d-4)/2})],
        \en
(b)
        \eq
        - \frac{\nabla^2 f_n(0;z_c)}{f_n(0;z_c)}
        =  v \sigma^2 n [1+{\cal O}(\beta n^{-\delta})],
        \en
(c)
        \eq
        \|\hat{D}^2 f_n(\cdot;z_c)\|_1
        \leq \mbox{const.} \beta n^{-d/2}.
        \en
The
error estimate in (a) is uniform in $\{k \in \Rd: a(k/\sqrt{v \sigma^2 n}) \leq
\gamma n^{-1} \log n
\}$, for $\gamma$ sufficiently small.
\end{theorem}

Part~(b) is an expression of diffusive behaviour.
Part~(c) is the first bound of \refeq{fbds}, and the factor
$n^{-d/2}$ is associated with the probability of return to the
origin of simple random walk after $n$ steps.

In the proof of Theorem~\ref{thm-1}, we will establish the
bounds of \refeq{fbds}
for all $m \in \Nbold$, with $z$ in an $m$-dependent interval containing
$z_c$.  Consequently, all bounds appearing
in Assumptions~E and $G$ follow as a Corollary, for $z=z_c$ and all $m\in
\N$.

We will also show that $z_c$, $A$ and $v$ obey the equations
    \eqarray
    1 & = &  \sum_{m = 1}^\infty g_m(0;z_c) , \label{z} \\
    A & = & \frac{1+ \sum_{m = 1}^\infty e_m(0;z_c)} {\sum_{m =
    1}^\infty m g_m(0;z_c)},   \label{A}\\
    v & = &-\frac{\sum_{m = 1}^\infty \nabla^2
    g_m(0;z_c)}{\sigma^2 \sum_{m = 1}^\infty m g_m(0;z_c)}.  \label{D}
    \enarray
It then follows immediately from the bounds of Assumptions~E and $G$
that
\eq
\lbeq{zAD}
    z_c=1+ \Ocal(\beta), \quad A=1+ \Ocal(\beta), \quad
    v = 1+ \Ocal(\beta).
\en
Equation (\ref{z}) states that at the critical value $z_c$ of the
recursion equation \refeq{fkrec}, with $k=0$, the coefficients
sum up to $1$ as in renewal theory. However, unlike the standard
situation in renewal theory, here the coefficients $g_m$ depend on
an additional variable $k$, and are not assumed to be
non-negative.
For the case of spread-out self-avoiding walk, the above bound for
$z_c$ improves on the error estimates of \cite{HHS01a,MS93,Penr94}.

With modest additional assumptions, the critical point $z_c$ can be
characterised in terms of
the {\em susceptibility}\/
\eq
\lbeq{sus}
    \chi(z) = \sum_{n=0}^\infty f_n(0;z).
\en

\begin{theorem}
\label{thm-zcpc}
Suppose that $d>4$ and that $L$ is sufficiently large.
Suppose there is a $p_c>0$ such that the susceptibility is finite and positive
for
$z \in (0,p_c)$, with $\lim_{z \uparrow p_c}\chi(z) = \infty$.
Suppose also
that the bounds of \refeq{fbds} for $z=z_c$ and all $m \geq 1$ imply the bounds of
Assumptions~E
and G for all $m \geq 2$, uniformly in $z \in [0,z_c]$.  Then $z_c=p_c$.
\end{theorem}

Next, we state a local central limit theorem.  Let
    \eq
    \lbeq{pndef}
    p_n(x;z) = \int_{[-\pi,\pi]^d} \frac{d^dk}{(2\pi)^d}
    e^{-ik\cdot x} f_n(k;z)
    \en
denote the inverse Fourier transform of $f_n(k;z)$.
The pointwise form
    \eq\lbeq{LCLT2}
    p_n(x\sqrt{v \sigma^2 n};z_c) \approx
    A\left(\frac{d}{2\pi n v \sigma^2 }\right)^{d/2}
    e^{-dx^2/2}
    \en
of the local central limit theorem
cannot hold for all solutions of the recursion \refeq{fkrec}.
For example, for self-avoiding walks starting
at the origin,
$p_n(0;z)=0$ for every $n \geq 1$
(see Section~\ref{sec-saw}).
This local
effect will disappear only if we average over a region that grows with
$n$.
For the averaging, we denote the cube of radius $R$ centred at $x \in \Zd$
by
\eq
\lbeq{CRdef}
    C_{R}(x) = \{y\in \Z^d: \|x-y\|_{\infty} \leq R\}.
\en
We use $\lfloor x\rfloor$ to denote
the closest lattice point in $\Z^d$ to $x\in \R^d$
(with an arbitrary rule to break ties), and write $a_n \sim b_n$
to mean $\lim_{n \to \infty }a_n/b_n =1$.

\begin{theorem}
\label{thm-3} Let $d>4$, and let $R_n$ be any
sequence with $\lim_{n\to \infty}
R_n = \infty$ and
$R_n = {o}(n^{1/2})$.
Then for all $x\in \R^d$ with $x^2[\log R_n]^{-1}$
sufficiently small,
    \eq
    \lbeq{LCLT}
    \frac{1}{(2R_n+1)^d }\sum_{y\in C_{R_n}(\lfloor x\sqrt{v \sigma^2 n}\rfloor)}
    p_n(y;z_c)
    \sim
    A\left(\frac{d}{2\pi n v \sigma^2 }\right)^{d/2} e^{-dx^2/2} .
    \en
\end{theorem}

An explicit error bound for \refeq{LCLT} in terms of
$R$ and $n$ is given in \refeq{LCLTerr}.
Note that \refeq{LCLT} is not a direct consequence of the
convergence of the Fourier transform in Theorem \ref{thm-1}. The
statement that sums over sets of volume $n^{d/2}$ converge to
integrals of the Gaussian density over the scaled set follows from
Theorem \ref{thm-1}. However, since we allow arbitrarily slow growth of
$R_n$ in Theorem~\ref{thm-3}, we are investigating $p_n$ on a smaller scale.

\subsection{Self-avoiding walks and oriented percolation}
\label{sec-examples}

In this section, we discuss two examples where our assumptions,
and hence our main results, apply.  The examples are spread-out
self-avoiding walks on $\Zd$, and critical spread-out oriented percolation on
$\Zd \times \Z_+$, in dimensions $d>4$.  Both models are based on a function
$D=D_L$ obeying the assumptions of Assumption~D.  Assumption~S is immediate
for both models.   We have two goals in this section.  The first goal
to provide a short sketch of the lace expansion, in order
to indicate how a recursion relation of the form \refeq{fkrec}
arises.  The second goal is to give some motivation for Assumptions~E
and G, as well as for our main results.

\subsubsection{Self-avoiding walks}
\label{sec-saw}

An $n$-step walk is a
mapping $\omega : \{0,\ldots, n\} \to \Zd$.
An $n$-step self-avoiding walk is an $n$-step walk with $\omega(i) \neq
\omega(j)$ for all $i \neq j$.  Let $\Rcal_n(x)$ denote the set of
$n$-step walks with $\omega(0)=0$ and $\omega(n)=x$,
and let ${\cal C}_n (x)$ denote
the subset of $\Rcal_n(x)$ consisting of
$n$-step self-avoiding walks.  For $\omega \in \Rcal_n(x)$, let
\eq
    W(\omega) = \prod_{i=1}^n D(\omega(i)-\omega(i-1)).
\en
We define $c_0(x) = \delta_{0,x}$, and, for $n \geq 1$,
\eq
\lbeq{cndef}
    c_n(x)=\sum_{\omega\in {\cal C}_n (x)} W(\omega).
\en
Without loss of generality, we now require that $D(0)=0$.

If the sum in \refeq{cndef}
were over $\Rcal_n(x)$, rather than $\Ccal_n(x)$, we would have
the recursion $c_{n+1}(x) = \sum_{y \in \Zd}c_1(y)c_n(x-y)$.
However, the right side of this equation includes contributions
from walks which visit the origin twice, and these are not present
in $c_{n+1}(x)$.
There should therefore be a correction term $F^{(1)}_{n+1}(x)$, defined by
the equation
\eq
\lbeq{ie.1}
    c_{n+1}(x) = \sum_{y \in \Zd}c_1(y)c_n(x-y) - F^{(1)}_{n+1}(x).
\en
The lace expansion is an expansion of $F^{(1)}_{n+1}(x)$
that treats self-avoiding walk as a small
perturbation of ordinary random walk for $d>4$.  It was
originally introduced in \cite{BS85} using a graph theoretic notion
of ``lace.'' Here, we follow instead the equivalent
inclusion-exclusion approach of \cite{Slad91} (see also \cite{MS93}),
for which the laces do not
explicitly appear.

To analyse the correction term $F^{(1)}_{n+1}(x)$, we define
$\Pcal_{n+1}(x)$ to be the set of $\omega \in \Rcal_{n+1}(x)$ for which
there
exists an $l \in \{2,\ldots,n+1\}$ (depending on $\omega$) with
$\omega(l)=0$ and $\omega(i) \neq \omega(j)$ for all
$i \neq j$ with $\{i,j\} \neq \{0,l\}$.
For the special case $x=0$, $\Pcal_{n+1}(0)$ is the set of $(n+1)$-step
self-avoiding polygons.  For general $x$, $\Pcal_{n+1}(x)$ is the set of
self-avoiding polygons followed by a self-avoiding walk from $0$ to $x$,
with the total length being $n+1$ and with the walk and polygon mutually
avoiding.
The set $\Pcal_{n+1}(x)$ is exactly the set of
walks that contribute to $\sum_{y \in \Zd}c_1(y)c_n(x-y)$
but do not contribute to $c_{n+1}(x)$.  Therefore
\eq
    F^{(1)}_{n+1}(x) = \sum_{\omega \in \Pcal_{n+1}(x)}W(x).
\en
Equation  \refeq{ie.1} can then be understood as just the
inclusion-exclusion relation: The first
term on the right side includes all walks from $0$ to $x$ which are
self-avoiding {\em after} the first step, and the second subtracts
the contribution due to those which are not self-avoiding from the
beginning, i.e., walks that return to the origin.


The inclusion-exclusion relation can now be applied to
$F^{(1)}_{n+1}(x)$.  We neglect the mutual avoidance of the polygon
and self-avoiding walk portions of an $\omega \in \Pcal_{n+1}(x)$,
and then correct by excluding the configurations which included intersections
of these two portions.
For $y \in \Zd$, let
\eq
    \pi_{m}^{(1)} (y) = \delta_{0,y}\sum_{\omega \in {\cal
    P}_m(0)}W(\omega).
\en
Then we define $F^{(2)}_{n+1}(x)$ by
\eq
    F^{(1)}_{n+1}(x)
    =
    \sum_{y \in \Zd} \sum_{m=2}^{n+1} \pi_m^{(1)}(y) c_{n+1-m}(x-y)
    - F^{(2)}_{n+1}(x).
\en
The correction term $F^{(2)}_{n+1}(x)$ involves configurations consisting
of a self-avoiding polygon and a self-avoiding walk from $0$ to $x$,
of total length $n+1$, and with an intersection required between the
self-avoiding polygon and the self-avoiding walk (in addition to their
intersection at the origin).
If we define $\Theta_{n+1}(x)$ to be the corresponding subset of
$\Rcal_{n+1}(x)$, then we have
\eq
    F^{(2)}_{n+1}(x) = \sum_{\omega \in \Theta_{n+1}(x)}W(\omega).
\en

The walk and polygon may intersect more than once, and we focus on the
first time an intersection occurs, measuring time according to the walk.
We then perform inclusion-exclusion again, neglecting the avoidance
between the portions of the self-avoiding walk before and after this
first intersection, and then subtracting a correction term.
The process is continued indefinitely.   We define
\eq
\lbeq{saw.pi2}
        \pi_{m}^{(2)} (y) =
        \sumtwo{m_1,m_2,m_3 \geq 1}{m_1+m_2+m_3=m}
        \prod_{j=1}^{3} \sum_{\omega_j \in {\cal C}_{m_j}(y)} W(\omega_i)
         I(\omega_1,\omega_2,\omega_3),
\en
where $I(\omega_1,\omega_2,\omega_3)$ is equal to $1$ if the
$\omega_i$ are pairwise mutually avoiding apart from their common
endpoints, and otherwise equals $0$.   Then we are led to
\eq
\lbeq{ckSAW}
    c_{n+1}(x) = \sum_{y \in \Zd}c_1(y)c_n(x-y) +
    \sum_{y \in \Zd} \sum_{m=2}^{n+1} \pi_m(y) c_{n+1-m}(x-y),
\en
with
        \eq
        \lbeq{lace/ie.pi}
        \pi_m (y) = \sum_{N=1}^{\infty}(-1)^{N} \pi_{m}^{(N)} (y),
        \en
with expressions for the $\pi^{(N)}$, for $N \geq 3$, which we do not
make explicit here.  The $\pi^{(N)}$ have convenient representations in terms
of Feynman diagrams \cite{BS85,MS93}.

By definition, $c_1(x)=D(x)$.  Using this in the first term of \refeq{ckSAW},
we then multiply through by $z^{n+1}$, with $z \geq 0$,
and take the Fourier
transform to obtain
        \eq
        \lbeq{ckSAWzeta}
        \hat{c}_{n+1}(k)z^{n+1} =
        z \hat{D}(k) \hat{c}_{n}(k) z^n
        + \sum_{m=2}^{n+1}
        \hat{\pi}_m(k)z^m \hat{c}_{n+1-m}(k)z^{n+1-m}.
        \en
Equation \refeq{ckSAWzeta} is a special case of \refeq{fkrec}, with
$e_n = 0$ and
\eq
\lbeq{zfegSAW}
    f_n(k;z) = \hat{c}_n(k) z^n \quad (n \geq 0),
    \quad
    g_1(k;z) = z \hat{D}(k), \quad
        g_n(k;z) =  \hat{\pi}_{n}(k)z^{n} \quad
        (n \geq 2).
\en
Assumption~S is a consequence of the lattice symmetry of
self-avoiding walks and the fact that $\hat{c}_n(k) \leq 1$.  We are
already assuming that $D$ obeys the requirements of Assumption~D, so
Assumption~D is immediate from \refeq{zfegSAW}.  Assumption~E is vacuous,
and the remaining assumption, Assumption~G, is the substantial one.

An analysis corresponding to the verification of some version of
Assumption~G is
present in one form or another in all lace expansion analyses
\cite{BS85,GI92,HS90a,HS90b,HS92b,HHS98,KLMS95,MS93,NY93,Slad87}.
In \cite{HS01b}, we
prove that Assumption~G is satisfied, as part of a larger
analysis proving Gaussian behaviour for networks of
mutually-avoiding spread-out self-avoiding walks in dimensions
$d>4$.

To give an indication how Assumption~G can be verified, we briefly discuss
this for the bound on $g_m$ arising from the contribution
due to $\hat{\pi}_m^{(2)}$.
For this, we simply neglect the mutual avoidance of the $\omega_j$
in \refeq{saw.pi2} to obtain an upper bound
\eq
\lbeq{saw.pi2bd}
    |\hat{\pi}_m^{(2)}(k)| \leq
    \sum_x \pi_m^{(2)}(x)
    \leq \sum_x \sumtwo{m_1,m_2,m_3 \geq 1}{m_1+m_2+m_3=m}
    \prod_{j=1}^{3} c_{m_j}(x).
\en
This succeeds in bounding $\hat{\pi}_m^{(2)}$ in terms of $c_l$ with $l<m$,
and the verification of Assumption~G can then be carried out, as explained
in detail in \cite{HS01b}.
The term
$\hat{\pi}^{(1)}_m$ is simpler, and higher order terms produce
powers of $\beta$ to ensure convergence of the sum over $N$ in
\refeq{lace/ie.pi}.  The restriction $d>4$ arises in bounding expressions
like the right side of \refeq{saw.pi2bd}, using \refeq{fbds}.

As a consequence of the verification of Assumption~G in \cite{HS01b},
Theorems~\ref{thm-1} and \ref{thm-3} hold for
$f_n(k;z_c) = \hat{c}_n(k) z_c^n$, with $z_c$, $A$
and $D$ given by Equations~(\ref{z})--(\ref{D}).
Theorem~\ref{thm-1}(a,b) were proved in \cite{MS93} with weaker error bounds,
via a different approach using generating functions.
Theorem~\ref{thm-3} is the first time a local central limit theorem
has been proved for the self-avoiding walk.
The limit $\lim_{n \to 0}\hat{c}_n(0)^{1/n}$ exists by a well-known
subadditivity argument \cite{MS93}.  Thus the critical value $z_c$ appearing
in Theorem~\ref{thm-1}(a) must be given by $z_c = \lim_{n \to
0}\hat{c}_n(0)^{-1/n}$, consistent with the conclusion of
Theorem~\ref{thm-zcpc}.
For uniform $D$, $z_c$ is thus the reciprocal
of the connective constant.

The first bound assumed in \refeq{fbds}, which asserts that
$\|\hat{D}^2 \hat{c}_m\|_1 \leq \Ocal(\beta z_c^{-m} m^{-d/2})$,
actually implies a bound on $\|c_m\|_\infty$. In
fact, if $m \geq 2$,
by ignoring the self-avoidance constraint for the first two
steps we have
        \eq
        \lbeq{bound-sup}
        \|c_{m}\|_{\infty}
        \leq \|D*D*c_{m-2}\|_{\infty}
        \leq \|\hat{D}^2 \hat{c}_{m-2}\|_1.
        \en
Since $z_c^{-1}\leq \Ocal(1)$ by \refeq{zAD}, the first bound of
\refeq{fbds} implies that
        \eq
        \lbeq{bound-sup1}
        \|c_m\|_\infty
        \leq \Ocal(\beta z_c^{-m} m^{-d/2}) .
        \en
The above argument gives \refeq{bound-sup1} when $m \geq 2$, but
\refeq{bound-sup1} clearly also holds for $m=1$ since
$\|c_1\|_\infty = \|D\|_\infty \leq \Ocal(\beta)$
and $z_c^{-1}$ is bounded away from zero.
The inequality \refeq{bound-sup1} is useful to estimate \refeq{saw.pi2bd}.
After having verified Assumption~G, the
bound $\|\hat{D}^2 \hat{c}_m\|_1 \leq \Ocal(\beta z_c^{-m} m^{-d/2})$
then follows from Theorem~\ref{thm-1}(c).  We conclude that
\refeq{bound-sup1}, which is of independent interest, does in fact
hold for all $m$.
An alternate proof of this bound, using generating functions and a finite
memory cut-off, is given in \cite{MS93} for the uniform $D$.

\subsubsection{Oriented percolation}

\newcommand{\nc}        { \conn  {\hspace{-2.8ex} /} \hspace{1.8ex}   }
\newcommand{\ndbc}      {\dbc {\hspace{-2.8ex} /} \hspace{+1.0ex} }

Our second example is spread-out oriented percolation.  The
setting is the graph with vertices $\Zd \times {\Zbold}_+$ and
directed bonds $((x,n),(y,n+1))$, for $n \geq 0$ and
$x,y \in \Zd$.   Let $z \in [0,\|D\|_\infty^{-1}]$, so that
$zD(y-x)\leq  1$.  We
associate to each directed bond $((x,n),(y,n+1))$
an independent random variable
taking the value $1$ with probability $zD(y-x)$ and $0$ with probability
$1-zD(y-x)$.
We say a bond is {\em occupied}\/ when the corresponding
random variable takes the value $1$, and {\em vacant}\/ when the
random variable is $0$.  The joint probability distribution of the
bond variables will be denoted ${\Pbold}_z$. We say that $(x,n)$
is {\em connected}\/ to $(y,m)$, and write $(x,n) \conn (y,m)$,
when there is an oriented path from $(x,n)$ to $(y,m)$ consisting
of occupied bonds, or if $(x,n)=(y,m)$. Given a bond configuration,
we define a bond to be {\em pivotal}\/ for $(x,n) \conn
(y,m)$ if $(x,n) \conn (y,m)$ when this bond is
made occupied, whereas $(x,n) \nc (y,m)$ when this bond is
made vacant. Finally, $(x,n)$ is {\em doubly connected}\/ to $(y,m)$,
denoted $(x,n) \dbc (y,m)$ when $x \conn y$ but
there is no pivotal bond for the connection from $(x,n)$
to $(y,m)$.  We write $(x,n) \ndbc (y,m)$ for the complement
of $(x,n) \dbc (y,m)$.

Given $z \in [0,\|D\|_\infty^{-1}]$, $n \geq 0$ and $x \in \Zd$,  we define the
two-point function
    \eq
    \lbeq{taukop}
    \tau_n(x;z) = {\Bbb P}_z((0,0) \conn (x,n)).
    \en
The lace expansion is an expansion for the two-point function.
It was derived specifically in the context of oriented percolation
by Nguyen and Yang in \cite{NY93,NY95},
using laces as in \cite{BS85}.  An expansion for
ordinary percolation based on inclusion-exclusion was derived in \cite{HS90a},
which applies also to oriented percolation and which is utilised in
\cite{HHS01b}.
Here, we briefly outline a derivation of the Nguyen--Yang
expansion using inclusion-exclusion.  This approach, which is different
than the inclusion-exclusion approach of \cite{HS90a}, is developed also in
\cite{Saka00}.

Given a configuration in which $(0,0) \conn (x,n)$, either $(0,0) \dbc (x,n)$
or $(0,0) \ndbc (x,n)$.  For $n \geq 0$, let
\eqarray
    \rho_n^{(0)}(x;z) &=& \Pbold_z((0,0) \dbc (x,n)), \\
    \sigma_n^{(0)}(x;z) &=& \Pbold_z(\{(0,0) \conn (x,n)\}
    \cap \{(0,0) \ndbc (x,n)\}) .
\enarray
Then
\eq
    \tau_{n+1}(x;z) =
        \rho_{n+1}^{(0)}(x;z) + \sigma_{n+1}^{(0)}(x;z) \quad (n \geq 0).
\en

A configuration contributing to $\sigma_{n+1}^{(0)}(x;z)$ contains
at least one pivotal bond, and hence a first pivotal bond $((u,m),(v,m+1))$.
Let $E$ denote the event that $(0,0) \dbc (u,m)$,
that $((u,m),(v,m+1))$ is occupied,
and that $(v,m+1) \conn (x,n+1)$.  Let $F$ be the event that the bond
$((u,m),(v,m+1))$ is pivotal for $(0,0) \conn (x,n+1)$.  Then
\eq
    \sigma_{n+1}^{(0)}(x;z) =
    \sum_{u,v} \sum_{m=0}^{n-1} \Pbold_z (E \cap F).
\en
We use inclusion-exclusion to write $\Pbold_z (E \cap F) = \Pbold_z (E )
- \Pbold_z (E \cap F^c)$.  Oriented percolation enjoys a Markov property
(unlike non-oriented percolation) which implies that
\eq
    \Pbold_z (E ) = \rho_m^{(0)}(u;z) \, zD(v-u) \, \tau_{n-m}(x-v;z).
\en
Thus we have
\eq
    \tau_{n+1}(x;z) =
        \rho_{n+1}^{(0)}(x;z) +
        \sum_{u,v} \sum_{m=0}^{n-1} \rho_m^{(0)}(u;z) \, zD(v-u) \,
        \tau_{n-m}(x-v;z)
        - \sigma_{n+1}^{(1)}(x;z),
\en
where
\eq
    \sigma_{n+1}^{(1)}(x;z)
        =\sum_{u,v} \sum_{m=0}^{n-1} \Pbold_z (E \cap F^c).
\en

For a configuration contributing to $E \cap F^c$, the first
pivotal bond for $(0,0) \conn (x,n+1)$, if there is one, will by definition
occur
later than time $m$.  We then perform inclusion-exclusion with respect to
the pivotal nature of this first pivotal bond, as above.  The procedure
is repeated indefinitely.  We pay special attention to the fact that
a site is always doubly connected to itself, so that
$\rho_0(x;z) = \delta_{0,x}$,
and we define $\pi^{(0)}_m(x;z) = \rho^{(0)}_m(x;z)$ for $ m\geq 2$,
and $\pi^{(0)}_m(x;z) = 0$ for $ m= 0,1$.
Repeating the inclusion-exclusion procedure outlined above then leads,
after taking the Fourier transform, to
    \eq \lbeq{taukrec}
    \hat{\tau}_{n+1}(k;z) = z\hat{D}(k) \tau_{n}(k;z)
        + z \hat{D}(k)\sum_{m=2}^{n}\hat{\pi}_m(k;z)\hat{\tau}_{n-m}(k;z)
        +\hat{\pi}_{n+1}(k;z),
    \quad \hat{\tau}_0(0;z)=1.
    \en
The function $\pi_m(x;z)$ is given by an explicit alternating series
\cite{HHS01b,NY93}
\eq
\lbeq{piopdef}
        \pi_m(x;z) = \sum_{N=0}^{\infty} (-1)^N \pi_m^{(N)}(x;z).
\en
Equation \refeq{taukrec} is a special case of \refeq{fkrec} with the choices
        \eq
        \lbeq{zfeOP}
        f_n(k;z)  = \hat{\tau}_{n}(k;z), \quad
        e_{n}(k;z) =\hat{\pi}_{n}(k;z),
        \en
and
        \eq
        \lbeq{gOP}
        g_1(k;z) = z \hat{D}(k), \quad g_2(k;z)=0, \quad
        g_m(k;z) = z\hat{D}(k)\hat{\pi}_{m-1}(k;z) \quad (m \geq 3).
        \en

For oriented percolation,
Assumption~S follows from the lattice symmetry of the model and
the fact that $\hat{\tau}_n(k;z) \leq (2n+1)^d$.  Assumption~D
follows from \refeq{zfeOP}--\refeq{gOP} and the fact that
$\hat{\pi}_1(k;z)=0$.  Assumptions~E and G are very closely related.
They can be verified, for $d>4$ and $L$ sufficiently large, by using
the BK inequality to bound $\pi^{(N)}$ in terms of the two-point
function itself.
For example, the BK inequality implies that $\pi_m^{(0)}(x;z) \leq
\tau_m(x;z)^2$ ($m \geq 2$).
This is carried out in detail in \cite{HHS01b}, where
the verification of the assumptions is an essential part of
a larger analysis relating the scaling limit of critical oriented
percolation for $d>4$ to super-Brownian motion.

Given the result of \cite{HHS01b} that
Assumptions~E and G do hold, the conclusions of
Theorems~\ref{thm-1} and \ref{thm-3} therefore apply at $z=z_c$.
Theorems~\ref{thm-1}(a,b) were proved in \cite{NY95}, with weaker error
bounds.  Theorems~\ref{thm-1}(c) and \ref{thm-3} are new, and the former
plays an essential role in the analysis of \cite{HHS01b}.
As was the case for self-avoiding walks in \refeq{bound-sup1}, it
is possible to conclude from the assumed bounds \refeq{fbds} that
    \eq
    \lbeq{tauinfbd}
    \| \tau_n (\cdot; z_c)\|_\infty \leq CK\beta n^{-d/2}.
    \en
See \cite{HHS01b} for details.

The hypotheses of Theorem~\ref{thm-zcpc} concerning the susceptibility
are well-known \cite{AN84}, and the hypothesis concerning Assumptions~E and G
is established in \cite{HHS01b}.  It therefore follows that
$z_c$ corresponds to the
critical oriented percolation threshold.  This can be understood directly
from the conclusion of Theorem~\ref{thm-1}(a), which asserts
in particular that $\lim_{n \to
\infty}\hat{\tau}_n(0;z_c) = A$, with $A$ positive and finite. The quantity
$\hat{\tau}_n(0;z)$ represents the expected
number of sites in the intersection of the connected cluster of
the origin with $\Zd \times \{n\}$.  This goes to zero for $z$ below the
percolation threshold,
by exponential decay of connectivities \cite{AB87,AN84,Mens86}.
On the other hand, for $z$ above the percolation threshold,
it can be expected that $\hat{\tau}_n(0;z) \to \infty$.
This latter statement is a consequence of the
shape theorem of \cite{BG90}.  The shape theorem has been proven only for
the nearest-neighbour model, but it can be expected to hold also
for the spread-out models.

The triangle condition \cite{AN84,BA91} is a diagrammatic
sufficient condition for the existence of several critical exponents.
This is discussed in detail in \cite{NY93}, where the triangle condition
is verified, in particular,
for sufficiently spread-out models of oriented percolation
when $d+1>5$.  The bound \refeq{tauinfbd}, together with the bound
$\sum_x \tau_n(x;z_c) < \infty$ (which is now proved, as discussed above),
can be combined to give an alternate proof of the triangle condition.
This is discussed in detail in \cite{HHS01b}.

\subsection{Organisation}

The remainder of this paper is organised as follows.  In
Section~\ref{sec-ih}, we introduce the induction hypotheses on
$f_n$ that will be used to prove our main results, and derive some
consequences of the induction hypotheses. The induction is advanced in
Section~\ref{sec-adv}. In Section~\ref{sec-pf}, the main results
stated in Section~\ref{sec-thm} are proved.  We conclude in
Appendix~\ref{sec-app}, by showing that the function
$D$ defined in \refeq{Dbasic2} obeys the hypotheses of Assumption~D.

\section{Induction hypotheses} \label{sec-ih}
\setcounter{equation}{0}

In this section, we introduce the induction hypotheses, verify that
they hold for $n=1$, discuss their motivation, and derive some of
their consequences.

\subsection{Statement of induction hypotheses (H1--H4)}
    \label{sec-ihstate}
The induction hypotheses involve a sequence $v_n$, which is
defined as follows.  We set $v_0=b_0=1$, and for $n \geq
1$ we define
    \eq
    b_n = -\frac{1}{\sigma^2}\sum_{m=1}^{n}
    \nabla^2 g_m(0;z)
    ,\quad
    c_n  =  \sum_{m=1}^{n} (m-1) g_m(0;z),\quad
    \label{Delta_n}
    v_n  = \frac{b_n}{1+c_n}.
    \en
The $z$--dependence of $b_n$, $c_n$, $v_n$ will usually be left implicit
in the notation.  We will often simplify the notation
by dropping $z$ also from $e_n$, $f_n$ and $g_n$, and write, e.g.,
$f_n(k)=f_n(k;z)$.

As we will explain in more detail in Section~\ref{sec-mot}, the diffusion
constant $\sigma^2 v$ of Theorem~\ref{thm-1} will turn out to be given by
$\sigma^2 v_\infty(z_c)$.  However, at this stage we have not yet
proved that the series in the definition of
$v_\infty$ converges.  Neither have we yet identified $z_c$, other than as a
solution to (\ref{z}), which involves a series whose convergence has not
yet been proved.

The induction hypotheses also involve several constants. Let $d>4$.  We fix
$\gamma, \delta, \rho >0$ according to
    \eq
    \lbeq{agddef}
    0< {\textstyle \frac{d-4}{2}} - \rho < \gamma
    < \gamma + \delta < 1 \wedge {\textstyle \frac{d-4}{2}}.
    \en
For example, we may first fix any $\delta \in (0, 1 \wedge
{\textstyle \frac{d-4}{2}})$ and then choose $\gamma$ and $\rho$
accordingly. We also introduce constants $K_1, \ldots , K_5$,
which are independent of $\beta$.

Let $z_0=z_1=1$, and define $z_n$ recursively by
    \eq
    \label{z_n}
    z_{n+1} = 1-\sum_{m=2}^{n+1}g_m(0;z_n),
    \qquad n \geq 1.
    \en
For $n \geq 1$, we define intervals
    \eq
    \lbeq{Indef}
    I_n = [z_n - K_1\beta n^{-(d-2)/2}, z_n + K_1\beta n^{-(d-2)/2}].
    \en
Recall the definition $a(k)=1-\hat{D}(k)$ from \refeq{adef}. Our
induction hypotheses are that the following four statements hold
for all $z \in I_n$ and all $1\leq j\leq n$.

\begin{description}
\item[(H1)]
$|z_j - z_{j-1}| \leq K_1 \beta j^{-d/2}$.
\item[(H2)]
$|v_j - v_{j-1}| \leq K_2 \beta j^{-(d-2)/2}$.
\item[(H3)]
For $k$ such that $a(k) \leq \gamma j^{-1}\log j$, $f_j(k;z)$ can
be written in the form
\[
    f_j(k;z) = \prod_{i=1}^j\left[
    1 -v_i a(k)+ r_i(k)  \right],
\]
with $r_i(k) = r_i(k;z)$ obeying
\[
    |r_i(0)|\leq K_3 \beta i^{-(d-2)/2},\quad
    |r_i(k)-r_i(0)| \leq K_3 \beta a(k) i^{-\delta}.
\]
\item[(H4)]
For $k$ such that $a(k) > \gamma j^{-1}\log j$, $f_j(k;z)$ obeys
the bounds
\[
    |f_j(k;z)| \leq K_4 a(k)^{-2-\rho}j^{-d/2}, \quad
    |f_j(k;z) -f_{j-1}(k;z)| \leq K_5 a(k)^{-1-\rho} j^{-d/2}.
\]
\end{description}

\medskip \noindent
Note that, for $k=0$, (H3) reduces to
$f_j(0) = \prod_{i=1}^j [1 + r_i(0) ]$.

We define
\eq
\lbeq{K3'def}
    K_4' = \max \{C_e(cK_4), C_g(cK_4), K_4\},
\en
where $c$ is a constant determined in Lemma~\ref{lem-pibds} below.
To advance the induction, we will need to assume that
    \eq
    \lbeq{Kcond}
    K_3 \gg K_1 > K_4' \geq K_4 \gg 1, \quad
    K_2 \geq K_1, 2K_4' , \quad
    K_5 \gg K_4.
    \en
Here $a \gg b$ denotes the statement that $a/b$ is sufficiently large.
The amount by which, for instance, $K_3$ must exceed $K_1$ is
independent of $\beta$ and will be determined during the course of
the advancement of the induction in Section~\ref{sec-adv}.

\subsection{Initialisation of the induction}

We now verify that the induction hypotheses hold when $n=1$.  Fix $z \in I_1$.

\begin{description}
\item[(H1)]
We simply have $z_1-z_0=1-1=0$.
\item[(H2)]
We simply have $|v_1-v_0| = |z-1|$, so that (H2) is satisfied provided
$K_2 \geq K_1$.
\item[(H3)]
We are restricted to $a(k)=0$.  By \refeq{Dbound1},
this means $k=0$. By
Assumption~D, $f_1(0;z) = z$, so that $r_1(0) =
z-z_1$.  Thus (H3) holds provided we
take $K_3 \geq K_1$.
\item[(H4)]
We note that $|f_1(k;z)| \leq z \leq 2$ for $\beta$
sufficiently small, $|f_1(k;z)-f_0(k;z)| \leq 3$, and $a(k) \leq
2$.  The bounds of (H4) therefore hold provided we take $K_4 \geq
2^{3+\rho}$ and $K_5 \geq 3 \cdot 2^{1+\rho}$.
\end{description}

\subsection{Discussion of induction hypotheses}
\label{sec-mot}

\noindent {\bf (H1) and the critical point.}\ The critical point
can be formally identified as follows. We set $k=0$ in \refeq{fkrec}, then
sum over $n$, and solve for the susceptibility $\chi(z)$
of \refeq{sus}.  The result is
    \eq
    \lbeq{chiz}
    \chi(z)
    = \frac{1+\sum_{m=2}^{\infty} e_m(0;z)}{1-\sum_{m=1}^\infty g_m(0;z)}.
    \en
The critical point should correspond to the smallest zero of the
denominator and hence should obey the equation
    \eq
    \lbeq{1pi}
        1 - \sum\limits_{m=1}^{\infty} g_m(0; z_c)
    = 1 - z_c - \sum\limits_{m=2}^{\infty} g_m(0; z_c) =0.
    \en
However, we do not know {\em a priori}\/ that the series in \refeq{chiz} or
\refeq{1pi} converge. We therefore approximate \refeq{1pi} with
the recursion (\ref{z_n}), which bypasses the convergence issue by
discarding the $g_m(0)$ for $m>n+1$ that cannot be handled at the
$n^{\rm th}$ stage of the induction argument. The sequence $z_n$
will ultimately converge to $z_c$. Equation~\refeq{1pi} is
identical to (\ref{z}).

In dealing with the sequence $z_n$, it is convenient to formulate
the induction hypotheses for a small interval $I_n$ approximating
$z_c$.  As we will see in Section \ref{sec-prel}, (H1) guarantees
that the intervals $I_j$ are decreasing: $I_1 \supset I_2 \supset
\cdots \supset I_n$. Because the length of these intervals is
shrinking to zero, their intersection $\cap_{j=1}^\infty I_j$ is a
single point, namely $z_c$. Hypothesis (H1) drives the convergence
of $z_n$ to $z_c$ and gives some control on the rate. The rate is
determined from (\ref{z_n}) and the ansatz that the difference
$z_j-z_{j-1}$ is approximately $-g_{j+1}(0,z_c)$, with
$|g_j(k;z_c)|= \Ocal( \beta j^{-d/2})$ as in  Assumption~G.

\medskip\noindent {\bf (H2).}\ The formula for
$v_n$ in (\ref{Delta_n}) can be motivated by the following rough
argument. Differentiating \refeq{fkrec} twice with respect to $k$,
setting $k=0$, and using the fact that odd derivatives vanish by
Assumption~S, we obtain
    \eq
    \lbeq{Alaprec}
    \nabla^2f_{n+1}(0) = \sum\limits_{m=1}^{n+1}
    \Big[g_m(0) \nabla^2 f_{n+1-m}(0)
    + \nabla^2 g_m(0) f_{n+1-m}(0) \Big]
    +\nabla^2e_{n+1}(0).
    \en
We will use `$\approx$' to denote an uncontrolled approximation
in a rough argument.
In \refeq{Alaprec}, we make the approximations $f_{n+1-m}(0)
\approx f_n(0)$ and $\nabla^2 e_{n+1} \approx 0$, subtract
$\nabla^2 f_n(0)$ from both sides, and recall the definition of
$b_n$ from (\ref{Delta_n}), to obtain
    \eq
    \lbeq{Alaprec2}
    \nabla^2f_{n+1}(0)  - \nabla^2 f_n(0)
    \approx  -\sigma^2 b_{n+1}f_n(0)
    + \sum\limits_{m=1}^{n+1}
    g_m(0) \nabla^2 f_{n+1-m}(0)
    -\nabla^2 f_n(0).
    \en
In view of (\ref{z_n}), we insert the approximation $1 \approx
\sum_{m=1}^{n+1} g_m(0)$ to obtain
    \eq
    \lbeq{Alaprec3}
    \nabla^2f_{n+1}(0) - \nabla^2 f_n(0)
    \approx -\sigma^2 b_{n+1}f_n(0)
    + \sum\limits_{m=1}^{n+1} g_m(0)
    \left[ \nabla^2 f_{n+1-m}(0) - \nabla^2 f_{n}(0) \right].
    \en

On the other hand, applying $\nabla^2$ to (H3) with $j=n$ and then
using $r_n(0) \approx 0$ and $\nabla^2 r_n(0) \approx 0$ gives
    \eq
    \lbeq{Dnmot}
    \nabla^2 f_n(0) - \nabla^2 f_{n-1}(0) \approx
    -\sigma^2 v_n f_{n-1}(0).
    \en
In view of this, we make the approximation $\nabla^2 f_{n+1-m}(0)-
\nabla^2 f_{n}(0) \approx (m-1)\sigma^2v_{n+1}f_n(0)$.  Recalling
the definition of $c_{n+1}$ from (\ref{Delta_n}), \refeq{Alaprec3}
then gives
    \eq
    \nabla^2f_{n+1}(0) - \nabla^2 f_n(0) \approx -\sigma^2
    (b_{n+1} - v_{n+1}c_{n+1})f_n(0).
    \en
Putting the right-hand side equal to $-\sigma^2v_{n+1}f_n(0)$, as
in \refeq{Dnmot}, leads to the formula for $v_{n+1}$ of
(\ref{Delta_n}).

The assumed bound on $|v_j - v_{j-1}|$ can then be guessed by
writing $v_j$ and $v_{j-1}$ in terms of the $b_i$ and $c_i$,
and assuming the bounds of Assumption~G to estimate the resulting
expression.  The calculations will be
carried out in detail when we advance (H2)
in Section~\ref{sec-advH1prime}.

\medskip\noindent {\bf (H3).}\
The bound on $r_i(0)$ ensures that the limit $A=\lim_{j \to
\infty}f_j(0;z_c) = \prod_{i=1}^\infty [1+r_i(0;z_c)]$ exists.
The expression for $f_j$ in (H3) can be approximately rewritten as
    \[
    f_j(k) \approx f_j(0) \exp{\left[\textstyle\sum_{i=1}^j
    \Big( -v_i a(k)+ r_i(k)-r_i(0) \Big) \right]} .
    \]
The bound on $r_i(k)-r_i(0)$ indicates
that it is a small perturbation of the leading term $-v_ia(k)$.
The limit $v = \lim_{n \to \infty}v_n$ exists, if we assume the
bounds on $g_m$ of Assumption~G.  Also, $a(k) \sim
\sigma^2 k^2/2d$ as $k \to 0$. Thus we can understand (H3) as a
precise version of the approximation
    \eq
    \lbeq{fjkapprox}
    f_j(k) \approx
    f_j(0)\exp \Big[- \frac{v\sigma^2 k^2 j}{2d}\Big] ,
    \en
which we expect to be valid at least for $\sigma^2 k^2 j$ of order
$1$. This is consistent with Theorem~\ref{thm-1}.
For large $j$, the restriction $a(k) \leq
\gamma j^{-1}\log j$ is essentially the restriction $\sigma^2 k^2
j \leq 2d \gamma \log j$, which includes the region where
$\sigma^2 k^2 j $ is order 1, plus some additional room to
manoeuvre.

\medskip\noindent {\bf (H4).}\
For $\sigma^2 k^2 j > 2d \gamma \log j$, we require (and can
prove) less accurate control of $f_j(k)$, as expressed in (H4).
The form of (H4) has been chosen in part to be less stringent than
(H3) for $a(k)=\gamma j^{-1}\log{j}$, where the transition from
(H3) to (H4) takes place. In fact, inserting $\sigma^2
k^2/2d = \gamma j^{-1}\log j$ into \refeq{fjkapprox} gives
an expression which grows like $j^{-\gamma v}$.  The
bounds on $g_m$ of Assumption~G would imply that $z_c =
1+\Ocal(\beta)$ and $v= 1+\Ocal(\beta)$, so that $j^{-\gamma v}=
j^{-\gamma(1+\Ocal(\beta))}$. On the other hand, putting
$a(k)=\gamma  j^{-1}\log j$ in the first bound of (H4) gives a
bound which grows like $j^{-(d-4)/2 + \rho}$ times a power of a
logarithm.  By \refeq{agddef}, this is a weaker bound than the
(H3) bound.

\subsection{Consequences of induction hypotheses}
\label{sec-prel} In this section we will derive some important
consequences of the induction hypotheses.
The key result is that the induction hypotheses imply \refeq{fbds}
for all $1 \leq m \leq n$, from which the bounds of Assumptions~E
and $G$ then follow, for $2 \leq m \leq n+1$.

Here, and throughout the
rest of this paper,
    \begin{itemize}
    \item[] $C$ denotes a strictly positive constant that may depend
    on $d,\gamma,\delta,\rho$, but {\it not}\/ on the $K_i$,
    {\it not}\/ on $k$, {\it not}\/ on $n$, and {\it not}\/ on
    $\beta$ (provided $\beta$ is sufficiently small, possibly
    depending on the $K_i$). The value of $C$ may change
    from line to line.
    \end{itemize}
The first lemma shows that the intervals $I_j$ are nested,
assuming (H1).

\begin{lemma}
\label{lem-In} Assume (H1) for $1 \leq j \leq n$. Then $I_1
\supset I_2 \supset \cdots \supset I_{n}$.
\end{lemma}

\proof
Suppose $z \in I_j$, with $2 \leq j \leq n$.  Then by
(H1) and \refeq{Indef},
    \eq
    |z-z_{j-1}| \leq |z-z_j| + |z_j - z_{j-1}|
    \leq \frac{K_1 \beta}{j^{(d-2)/2}}
    + \frac{K_1 \beta}{j^{d/2}}
    \leq \frac{K_1 \beta}{(j-1)^{(d-2)/2}},
    \en
and hence $z \in I_{j-1}$.
\qed

By Lemma~\ref{lem-In}, if $z \in I_j$ for $1 \leq j \leq n$, then
$z \in I_1$ and hence, by \refeq{Indef},
    \eq
    \lbeq{znear1}
    |z -1| \leq K_1 \beta .
    \en
It also follows from (H2)
that, for $z \in I_n$ and $1 \leq j \leq n$,
    \eq
    \lbeq{vnear1}
    |v_j -1| \leq CK_2 \beta.
    \en

The induction hypothesis (H3) has the useful alternate form
\eq
\lbeq{fs}
    f_j(k) = f_j(0)\prod_{i=1}^j
    \left[
    1 - v_i a(k) + s_i(k)
    \right],
\en
with $s_i(0)=0$ and
\eq
\lbeq{sbd}
    |s_i(k)| \leq K_3 (1+C(K_2+K_3)\beta) \beta a(k) i^{-\delta}.
\en
In fact, \refeq{fs} is an identity if we define
\eq
\lbeq{sdef}
    s_i(k) = [1+r_i(0)]^{-1} [v_i a(k) r_i(0) + (r_i(k) - r_i(0))],
\en
and \refeq{sbd} then follows from \refeq{vnear1}
and the bounds on $r_i$ of (H3).

The next lemma provides an important upper bound on $f_j(k;z)$, for $k$ small
depending on $j$, as in (H3).

\begin{lemma}
\label{lem-cA} Let $z\in I_n$ and assume (H2--H3) for $1 \leq
j \leq n$. Then for $k$ with $a(k) \leq \gamma j^{-1}\log j$,
\eq
    |f_j(k;z)| \leq e^{CK_3\beta} e^{-(1-C(K_2+K_3)\beta)ja(k)}.
\en
\end{lemma}

\proof
We use \refeq{fs}, and conclude from the bound on $r_i(0)$ of (H3)
that $|f_j(0)| \leq e^{CK_3\beta}$.  Then we use \refeq{vnear1} and
\refeq{sbd} to obtain
\eq
    \prod_{i=1}^j  \left| 1 - v_i a(k) + s_i(k) \right|
    \leq \prod_{i=1}^j  \left| 1 - (1-CK_2\beta) a(k) +
    CK_3 \beta a(k) i^{-\delta} \right|.
\en
The desired bound then follows, once we use $1+x \leq e^x$ for each
factor on the right side.
\qed

The middle bound of \refeq{fbds} follows, for $ 1 \leq m \leq n$
and $z \in I_m$, directly
from Lemma~\ref{lem-cA}.  We next prove
two lemmas which provide the other two bounds of \refeq{fbds}.
This will supply the hypothesis \refeq{fbds} for Assumptions~E and G,
and therefore plays a crucial role in advancing the induction.

\begin{lemma}
\label{lem-L1norm}
Let $z \in I_n$ and assume (H2), (H3) and (H4).
Then for $1 \leq j \leq n$,
\eq
    \| \hat{D}^2 f_j(\cdot ;z)\|_1 \leq C(1+K_4) \beta j^{-d/2}.
\en
\end{lemma}

\proof
Fix $z \in I_n$ and $1 \leq j \leq n$, and define
    \eqarray
    R_{1} & = &\{k \in
    [-\pi,\pi]^{d}: a(k) \leq \gamma j^{-1}\log j , \;
    \|k\|_\infty \leq L^{-1} \} ,\nonumber \\
    R_{2} & = &\{k \in
    [-\pi,\pi]^{d}: a(k) \leq \gamma j^{-1}\log j , \;
    \|k\|_\infty > L^{-1} \} ,\nonumber \\
    \lbeq{Rjdef}
    R_{3} & = &\{k \in
    [-\pi,\pi]^{d}: a(k) > \gamma j^{-1}\log j, \;
    \|k\|_\infty \leq L^{-1} \},\nn\\
    R_4 & = &\{k \in
    [-\pi,\pi]^{d}: a(k) > \gamma j^{-1}\log j, \;
    \|k\|_\infty > L^{-1} \}.
    \enarray
The set $R_2$ is empty if $j$ is sufficiently large.  Then
    \eq
    \|\hat{D}^2f_{j}\|_1 = \sum_{i=1}^4\int_{R_i} \hat{D}(k)^2|f_{j}(k)|
    \frac{d^dk}{(2\pi)^d}.
    \en
We will treat each of the four terms on the right side separately.

On $R_1$, we use \refeq{Dbound1} in conjunction with
Lemma~\ref{lem-cA} and the fact that $\hat{D}^2
\leq 1$, to obtain
    \eq
    \int_{R_1} \hat{D}(k)^2|f_{j}(k)| \frac{d^dk}{(2\pi)^d}
    \leq \int_{R_1} Ce^{-cj(Lk)^2} \frac{d^dk}{(2\pi)^d}
    \leq \frac{C}{L^d j^{d/2}}.
\en
On $R_2$, we use Lemma~\ref{lem-cA} and \refeq{Dbound2} to
conclude that there is an $\alpha > 1$ such that
    \eq
    \int_{R_2} \hat{D}(k)^2|f_{j}(k)| \frac{d^dk}{(2\pi)^d}
    \leq \int_{R_2} \alpha^{-j} \frac{d^dk}{(2\pi)^d}
    = \alpha^{-j} |R_2|,
\en
where $|R_2|$ denotes the volume of $R_2$.  This volume is maximal
when $j=3$, so that
    \eq
    |R_2| \leq |\{ k : a(k) \leq 1 - \textstyle\frac{\gamma \log 3}{3} \}|
    \leq |\{k: \hat{D}(k) \geq \textstyle\frac{\gamma \log 3}{3} \}|
    \leq (\textstyle\frac{3}{\gamma \log 3})^2 \|\hat{D}^2\|_1
    \leq (\textstyle\frac{3}{\gamma \log 3})^2 \beta.
\en
Therefore $\alpha^{-j}|R_2| \leq C\beta j^{-d/2}$ and
    \eq
    \int_{R_2} \hat{D}(k)^2|f_{j}(k)| \frac{d^dk}{(2\pi)^d}
    \leq C\beta j^{-d/2}.
    \en

On $R_3$ and $R_4$, we use (H4).  As a result, the contribution
from these two regions is bounded above by
    \eq
    \frac{K_4 }{j^{d/2}} \sum_{i=3}^4\int_{R_i}
        \frac{\hat{D}(k)^2}{a(k)^{2+\rho}} \frac{d^dk}{(2\pi)^d}.
    \en
On $R_3$, we use $\hat{D}(k)^2\leq 1$ and \refeq{Dbound1} to
obtain the upper bound
    \eq
    \label{largekext}
    \frac{CK_4}{j^{d/2} L^{4+2\rho}}
    \int_{\|k\|_\infty< L^{-1}} \frac{1}{k^{4+2\rho}} d^dk
    \leq
    \frac{CK_4}{j^{d/2} L^{4+2\rho}}
    \left(\frac{1}{L}\right)^{d-4-2\rho}
    = \frac{CK_4\beta}{j^{d/2}},
    \en
where the integral is finite since $\rho < \frac{d-4}{2}$ by
\refeq{agddef}. On $R_4$, we
use \refeq{Dbound2} to obtain the bound
    \eq
    \frac{CK_4}{j^{d/2}}
    \int_{[-\pi,\pi]^d}  \hat{D}(k)^2\frac{d^dk}{(2\pi)^d}
    \leq \frac{CK_4\beta}{j^{d/2}}.
    \en
This completes the proof.
\qed

\begin{lemma}
\label{lem-fder}
Let $z \in I_n$ and assume (H2) and (H3).  Then, for
$1 \leq j \leq n$,
\eq
    | \nabla^2 f_j(0 ;z) | \leq (1+(K_2 + K_3) \beta ) \sigma^2 j.
\en
\end{lemma}

\proof
Fix $z\in I_n$ and $j$ with $1 \leq j \leq n$.
By \refeq{fs} and Assumption~S,
    \eq
    \lbeq{1.2b2}
    \nabla^2 f_{j}(0)
    = f_j(0)\sum_{i=1}^{j} \bigl[-\sigma^2v_i+\nabla^2 s_i(0)\bigr].
    \en
By \refeq{vnear1}, $|v_i -1| \leq CK_2\beta$.  For the second term on the
right side, we let
$e_1,\ldots,e_d$ denote the standard basis vectors in
$\Rd$. By \refeq{sbd}, for all $i\leq n$ we have
    \eq
    \lbeq{1.2b3}
    |\nabla^2 s_i(0)| = \Big| \sum_{l=1}^d \lim_{t \rightarrow 0}
    \frac{s_i(te_l)-s_i(0)}{t^2} \Big|
    \leq CK_3 \beta i^{-\delta}
    \sum_{l=1}^d \lim_{t \rightarrow 0} \frac{a(te_l)}{t^2}
    =  CK_3 \sigma^2 \beta i^{-\delta}.
    \en
Therefore, by Lemma~\ref{lem-cA},
    \eq
    |\nabla^2 f_{j}(0)| \leq e^{CK_3\beta} \sigma^2 j
    \Big(1+C(K_2 + K_3 )\beta\Big).
    \en
This completes the proof.
\qed

The next lemma is the key to advancing the induction, as it
provides bounds for $e_{n+1}$ and $g_{n+1}$.

\begin{lemma}
\label{lem-pibds}
 Let $z\in I_{n}$, and assume (H2),
(H3) and (H4). For $k \in [-\pi,\pi]^d$, $2 \leq j \leq n+1$, and
$\epsilon' \in [0,\epsilon]$, the following hold:
\vspace{-2mm}
\begin{tabbing}
(iii) \= \kill (i) \>  $|g_j(k;z)|\leq  K_4' \beta j^{-d/2}$,
\\
(ii) \> $|\nabla^2 g_j(0;z)|\leq   K_4'  \sigma^2 \beta
j^{-(d-2)/2}$,
\\
(iii) \> $|\partial_z g_j(0;z)|\leq  K_4' \beta j^{-(d-2)/2},$
\\
(iv) \>  $|g_j(k;z)-g_j(0;z)- a(k) \sigma^{-2}\nabla^2 g_j(0;z)|
\leq  K_4' \beta a(k)^{1+\epsilon'}j^{-(d-2-2\epsilon')/2},$
\\
(v)   \> $|e_j(k;z)|\leq  K_4' \beta j^{-d/2}$,
\\
(vi)  \> $|e_j(k;z)-e_j(0;z)|\leq  K_4' a(k) \beta
j^{-(d-2)/2}.$
\end{tabbing}
\end{lemma}

\proof The bounds \refeq{fbds} for $1 \leq m \leq n$
follow from Lemmas~\ref{lem-cA}--\ref{lem-fder}, with
$K=cK_4$ (this defines $c$), assuming that $\beta$ is sufficiently
small. The bounds of the lemma then follow immediately from
Assumptions~E and G, with $K_4'$ given in \refeq{K3'def}.
\qed

\section{The induction advanced}
\label{sec-adv} \setcounter{equation}{0}
In this section we advance the induction hypotheses (H1--H4) from
$n$ to $n+1$.
Throughout this section, in accordance with the uniformity
condition on (H2--H4), we fix $z \in I_{n+1}$.  We also assume that
$\beta \ll 1$.

\subsection{Advancement of (H1)}
\label{sec-advH1}

By (\ref{z_n}) and the mean-value theorem,
    \eqarray
    z_{n+1}-z_n & = & -\sum_{m=2}^{n}
    [g_m(0;z_n)-g_m(0;z_{n-1})] -g_{n+1}(0;z_n)\nonumber \\
            & = &
    -(z_n-z_{n-1})\sum_{m=2}^{n} \partial_z g_m(0;y_n) -
    g_{n+1}(0;z_{n}),
    \enarray
for some $y_n$ between $z_n$ and $z_{n-1}$.  By (H1) and
\refeq{Indef}, $y_n \in I_n$. Using Lemma~\ref{lem-pibds} and
(H1), it then follows that
    \eqarray
    |z_{n+1}-z_n| &\leq& K_1\beta n^{-d/2}
    \sum\limits_{m=2}^{n}  K_4' \beta m^{-(d-2)/2} +  K_4' \beta
    (n+1)^{-d/2} \nonumber \\
    & &\nonumber \\
    &\leq& K_4' \beta(1 + CK_1 \beta) (n+1)^{-d/2}. \label{zdiff}
    \enarray
Thus (H1) holds for $n+1$, for $\beta$ small and $K_1 > K_4'$.

\medskip
Having advanced (H1) to $n+1$, it then follows from
Lemma~\ref{lem-In} that $I_1 \supset I_2 \supset \cdots \supset
I_{n+1}$.

For $n \geq0$, define
    \eq
    \lbeq{zetadef}
        \zeta_{n+1} =\zeta_{n+1}(z) = \sum_{m=1}^{n+1} g_m(0;z)-1
    = \sum_{m=2}^{n+1} g_m(0;z) + z -1.
    \en
The following lemma, whose proof makes use of (H1) for $n+1$, will
be needed in what follows.

\begin{lemma} \label{zetan} For all $z\in I_{n+1}$,
    \eqarray
    \lbeq{zetanbd}
    |\zeta_{n+1}| & \leq & CK_1 \beta (n+1)^{-(d-2)/2}.
    \lbeq{Op1}
    \enarray
\end{lemma}

\proof By (\ref{z_n}) and the mean-value theorem,
    \eqarray
    |\zeta_{n+1}| & = & \Big|(z-z_{n+1}) + \sum_{m=2}^{n+1}
    [g_m(0;z)-g_m(0;z_n)]\Big|\nonumber\\
              & = & \Big|(z-z_{n+1})
    + (z-z_n) \sum_{m=2}^{n+1} \partial_z g_m(0;y_n) \Big|,
    \enarray
for some $y_n$ between $z$ and $z_n$. Since $z\in I_{n+1} \subset
I_n$ and $z_n \in I_n$, we have $y_n \in I_n$. Therefore, by
Lemma~\ref{lem-pibds},
\eq
        |\zeta_{n+1}|  \leq K_1\beta(n+1)^{-(d-2)/2}
    + K_1\beta n^{-(d-2)/2}
    \sum_{m=2}^{n+1}  K_4' \beta m^{-(d-2)/2}
    \leq  K_1\beta(1+CK_4'\beta) (n+1)^{-(d-2)/2}.
\en
The lemma then follows, for $\beta$ sufficiently small.
\qed

\subsection{Advancement of (H2)}
\label{sec-advH1prime}

Let $z\in I_{n+1}$.  As observed in Section~\ref{sec-advH1}, this
implies that $z \in I_j$ for all $j \leq n+1$.
The definitions in (\ref{Delta_n}) imply that
    \eq
    \lbeq{vinc}
    v_{n+1} - v_n = \frac{1}{1+c_{n+1}}(b_{n+1}-b_n) -
    \frac{b_n}{(1+c_n)(1+c_{n+1})}(c_{n+1}-c_n),
    \en
with
    \eq
    \lbeq{bcdiff}
    b_{n+1}-b_n = -\frac{1}{\sigma^2}\nabla^2 g_{n+1}(0)
    , \quad
    c_{n+1}-c_n =  n g_{n+1}(0) .
    \en
By Lemma~\ref{lem-pibds}, both differences in \refeq{bcdiff} are bounded by
$K_4' \beta (n+1)^{-(d-2)/2}$, and, in addition,
\eq
\lbeq{bnear1}
    |b_j - 1 | \leq C K_4' \beta, \quad
    |c_j | \leq C K_4' \beta
\en
for $1 \leq j \leq n+1$.  Therefore
    \eq
    |v_{n+1} - v_n|\leq K_2 \beta (n+1)^{-(d-2)/2},
    \en
provided we assume $K_2 > 2K_4'$.  This advances (H2).

\subsection{Advancement of (H3)}
\label{sec-advH2}

\subsubsection{The decomposition}

The advancement of the induction hypotheses (H3--H4)
is the most technical part of the proof.
For (H3), we fix $k$ with
$a(k)\leq \gamma (n+1)^{-1}\log{(n+1)}$, and $z\in I_{n+1}$. The
induction step will be achieved as soon as we are able to write
the ratio $f_{n+1}(k)/f_n(k)$ as
\eq
    \frac{f_{n+1}(k)}{f_n(k)} = 1-v_{n+1}a(k) + r_{n+1}(k),
\en
with $r_{n+1}(0)$ and $r_{n+1}(k)-r_{n+1}(0)$
satisfying the bounds required by (H3).

To begin, we divide the recursion relation \refeq{fkrec} by
$f_n(k)$, and use \refeq{zetadef}, to obtain
    \eqarray
    \label{rec hat{tau}_n(k)}
    \frac{f_{n+1}(k)}{f_{n}(k)} & =
    &1 +\sum_{m=1}^{n+1} \Big[g_m(k)
    \frac{f_{n+1-m}(k)}{f_{n}(k)}-g_m(0)\Big] +\zeta_{n+1}+
    \frac{e_{n+1}(k)}{f_n(k)}.
    \enarray
By (\ref{Delta_n}),
\eq
    v_{n+1} = b_{n+1}-v_{n+1}c_{n+1}
    =-\sigma^{-2} \sum_{m=1}^{n+1}\nabla^2 g_m(0) - v_{n+1}
    \sum_{m=1}^{n+1}(m-1) g_m(0).
    \en
Thus we can rewrite (\ref{rec hat{tau}_n(k)}) as
    \eq
    \frac{f_{n+1}(k)}{f_{n}(k)} = 1 - v_{n+1} a(k) + r_{n+1}(k),
    \label{the eq}
    \en
where
    \eq
    r_{n+1}(k) = X(k)+Y(k)+Z(k)+ \zeta_{n+1}
    \en
with
    \eqarray
    X(k)    & = &\sum_{m=2}^{n+1}
    \Big[\big(g_m(k)-g_m(0)\big)\frac{f_{n+1-m}(k)}{f_n(k)}-a(k)
    \sigma^{-2}\nabla^2g_m(0)\Big],\\
    \lbeq{H2IIdef}
    Y(k) & = &\sum_{m=2}^{n+1} g_m(0)\left[
    \frac{f_{n+1-m}(k)}{f_n(k)}-1- (m-1) v_{n+1}a(k)\right],
    \hspace{6mm} \\
    Z(k) & = & \frac{e_{n+1}(k)}{f_n(k)}.
    \enarray
The $m=1$ terms in $X$ and $Y$ vanish and have not been
included.

We will prove that
\eq
\lbeq{rbds}
    |r_{n+1}(0)|\leq \frac{C(K_1+ K_4') \beta}{(n+1)^{(d-2)/2}},
    \qquad \qquad
    |r_{n+1}(k)-r_{n+1}(0)|
    \leq \frac{C K_4' \beta a(k)}{(n+1)^{\delta}}.
\en
This gives (H3) for $n+1$, provided we assume that $K_3 \gg K_1$
and $K_3 \gg K_4'$.
To prove the bounds on $r_{n+1}$ of \refeq{rbds}, it will be convenient
to make use of some elementary convolution bounds, as well as some bounds
on ratios involving $f_j$.  These preliminary bounds are given
in Section~\ref{sec-ratiobds}, before we present the proof of
\refeq{rbds} in Section~\ref{sec-XYZ}.

\subsubsection{Convolution and ratio bounds}
\label{sec-ratiobds}

The proof of \refeq{rbds} will make use of the following
elementary convolution bounds.  To keep the discussion simple, we
do not obtain optimal bounds.

\begin{lemma}
    \label{lem-conv}
For $n \geq 2$,
    \eq
    \lbeq{conv-bound}
    \sum_{m=2}^{n} \frac{1}{m^a}
    \sum_{j=n-m+1}^n \frac{1}{j^b} \leq \left\{
    \begin{array}{lll}&C n^{-a\wedge b+1} &\mbox{for }a,b>1\\
    &Cn^{-(a-2)\wedge b}    &\mbox{for }a>2, b>0\\
    &Cn^{-(a-1)\wedge b}    &\mbox{for }a>2, b>1\\
    &Cn^{-a\wedge b}    &\mbox{for }a,b>2.\end{array}\right.
    \en
\end{lemma}

\proof
Since $m+j\geq n$, either $m$ or $j$ is at least $\frac{n}{2}$.
Therefore
    \eq
    \lbeq{conv-bound2}
    \sum_{m=2}^{n} \frac{1}{m^a}
    \sum_{j=n-m+1}^n \frac{1}{j^b} \leq \left(\frac{2}{n}\right)^a
    \sum_{m=2}^{n}\sum_{j=n-m+1}^n \frac{1}{j^b}
    +\left(\frac{2}{n}\right)^b \sum_{m=2}^{n}\sum_{j=n-m+1}^n
    \frac{1}{m^a}.
    \en
If $a,b>1$, then the first term is bounded by $Cn^{1-a}$ and the
second by $Cn^{1-b}$.\\ If $a>2, b>0$, then the first term is
bounded by $Cn^{2-a}$ and the second by $Cn^{-b}$.\\ If $a>2,
b>1$, then the first term is bounded by $Cn^{1-a}$ and the second
by $Cn^{-b}$.\\ If $a,b>2$, then the first term is bounded by
$Cn^{-a}$ and the second by $Cn^{-b}$.
\qed

We also will make use of several estimates involving ratios.
We begin with some preparation.
Given a vector $x=(x_l)$ with $\sup_l|x_l|<1$, define $\chi(x) = \sum_l
\frac{|x_l|}{1-|x_l|}$.
The bound $(1-t)^{-1} \leq \exp [t(1-t)^{-1}]$, together with Taylor's
Theorem applied to $f(t) = \prod_l \frac{1}{1-tx_l}$, gives
\eq
\lbeq{Taylor1}
    \left|\prod_l \frac{1}{1-x_l} -1\right| \leq \chi(x) e^{\chi(x)},
    \quad
    \left|\prod_l \frac{1}{1-x_l} -1 -\sum_l x_l\right|
    \leq \frac 32\chi(x)^2 e^{\chi(x)}.
\en

We assume throughout the rest of this section that $a(k) \leq \gamma
(n+1)^{-1} \log (n+1)$ and $2 \leq m \leq n+1$, and define
\eq
    \psi_{m,n} = \sum_{j=n+2-m}^n \frac{|r_j(0)|}{1-|r_j(0)|}, \quad
\lbeq{chidef}
    \chi_{m,n}(k) = \sum_{j=n+2-m}^n \
    \frac{v_j a(k)+ |s_j(k)|}{1-v_j a(k) - |s_j(k)|}.
\en
By \refeq{vnear1} and \refeq{sbd},
    \eq
    \lbeq{chibd1}
    \chi_{m,n}(k) \leq (m-1) a(k)Q(k)
    \quad \mbox{ with } \quad Q(k) = [1+C(K_2+K_3)\beta][1+Ca(k)].
    \en
Since $a(k) \leq \gamma (n+1)^{-1} \log (n+1)
\leq \gamma (m-1)^{-1}\log (m-1)$, we have
$m-1 \leq \gamma a(k)^{-1} \log (m-1)$.  Therefore
    \eq
    \lbeq{chibd2}
    e^{\chi_{m,n}(k)} \leq e^{(m-1) a(k)Q(k)}
    \leq (m-1)^{\gamma Q(k)}
    = (m-1)^{\gamma Q(0)}(m-1)^{\gamma (Q(k)-Q(0))}.
    \en
Now $Q(k)-Q(0) \leq Ca(k)$, and hence, for any $q < Q(0)$ we have
\eq
\lbeq{chibd3}
    e^{\chi_{m,n}(k)}
    \leq (m-1)^{\gamma Q(0)} C a(k)^{-Ca(k)} (\log (m-1))^C
    \leq C (n+1)^{\gamma q}.
\en
Here, we have used the fact that $a(k)^{-Ca(k)}$ remains bounded as
$k \to 0$, and we have absorbed the logarithmic factor by a small power.
We can choose $q$ to be as close to $1$ as desired, by taking $\beta$ to
be small.

We now turn to the ratio bounds.  It follows from (H3)
and the first inequality of \refeq{Taylor1} that
\eq
\lbeq{ratio1.a}
    \left| \frac{f_{n+1-m}(0)}{f_n(0)} - 1 \right|
    \leq \psi_{m,n}e^{\psi_{m,n}}
    \leq
    \sum_{j=n+2-m}^m \frac{CK_3 \beta}{j^{(d-2)/2}}
    \leq
    \frac{CK_3 \beta}{(n+2-m)^{(d-4)/2}}.
\en
Therefore
\eq
\lbeq{ratio0}
    \left| \frac{f_{n+1-m}(0)}{f_n(0)} \right|
    \leq
    1+CK_3 \beta.
\en
By \refeq{fs},
\eq
    \left| \frac{f_{n+1-m}(k)}{f_n(k)} - 1 \right|
    \leq
    \left| \frac{f_{n+1-m}(0)}{f_n(0)} \right|
    \left| \prod_{j=n+2-m}^n [1-v_ja(k) +s_j(k)]^{-1}
    - 1 \right|
    + \left| \frac{f_{n+1-m}(0)}{f_n(0)} - 1 \right|.
\en
The first inequality of \refeq{Taylor1}, together with
\refeq{chibd1}--\refeq{ratio0}, then gives
\eq
\lbeq{ratio1}
    \left| \frac{f_{n+1-m}(k)}{f_n(k)} - 1 \right|
    \leq
    C (m-1) a(k)(n+1)^{\gamma q} + \frac{CK_3 \beta}{(n+2-m)^{(d-4)/2}} .
\en
Similarly,
\eq
\lbeq{ratio2}
    \left| \frac{f_{n}(0)}{f_n(k)} - 1 \right|
    \leq
    C a(k) (n+1)^{1+\gamma q}.
\en

Next, we estimate the quantity $R_{m,n}(k)$, which is defined by
\eq
\lbeq{Rdef}
    R_{m,n}(k) = \prod_{j=n+2-m}^n [1-v_ja(k) +s_j(k)]^{-1}
    - 1
    - \sum_{j=n+2-m}^n [v_ja(k) -s_j(k)] .
\en
By the second inequality of \refeq{Taylor1}, together with
\refeq{chibd1} and \refeq{chibd3}, this obeys
\eq
\lbeq{Rbd}
    |R_{m,n}(k)| \leq \frac{3}{2} \chi_{m,n}(k)^2 e^{\chi_{m,n}(k)}
    \leq
    C m^2 a(k)^2 (n+1)^{\gamma q}.
\en

Finally, we apply (H3) to obtain
\eq
\lbeq{ratio3}
    \left| \frac{f_{m-1}(k)}{f_m(k)} - 1 \right|
    =
    \left| [1-v_ma(k) +(r_m(k)-r_m(0)) + r_m(0)]^{-1}
    - 1 \right|
    \leq
    Ca(k) + \frac{CK_3 \beta}{m^{(d-2)/2}} .
\en

\subsubsection{The induction step}
\label{sec-XYZ}

By definition,
    \eq
    \lbeq{rp0}
    r_{n+1}(0) = Y(0)+Z(0)+\zeta_{n+1}
    \en
and
    \eq
    \lbeq{rpk0}
    r_{n+1}(k)-r_{n+1}(0) = X(k) + \Big( Y(k)-Y(0) \Big)
     + \Big( Z(k)-Z(0) \Big).
    \en
Since $|\zeta_{n+1}| \leq CK_1 \beta (n+1)^{-(d-2)/2}$ by
Lemma~\ref{zetan}, to prove \refeq{rbds} it suffices to show that
    \eq
    \lbeq{rp0suf}
    |Y(0)| \leq C K_4' \beta (n+1)^{-(d-2)/2} ,\quad
    |Z(0)| \leq C K_4' \beta (n+1)^{-(d-2)/2}
    \en
and
    \eqarray
    \lbeq{rpk0suf}
    && |X(k)|  \leq  C K_4' \beta a(k) (n+1)^{-\delta}, \quad
    |Y(k)-Y(0)| \leq C K_4' \beta a(k) (n+1)^{-\delta},
    \nonumber \\ && \hspace{25mm}
    |Z(k)-Z(0)|  \leq  C K_4' \beta a(k) (n+1)^{-\delta}.
    \enarray
The remainder of the proof is devoted to establishing
\refeq{rp0suf} and \refeq{rpk0suf}.

\medskip\noindent
{\em Bound on $X$}. We write $X$ as $X=X_1+X_2$, with
    \eqarray
    \lbeq{X1def}
    X_1 &=& \sum_{m=2}^{n+1}
    \Big[g_m(k)-g_m(0)-a(k) \sigma^{-2}\nabla^2g_m(0)\Big] ,
    \\
    X_2 &=&\sum_{m=2}^{n+1}
    \Big[g_m(k)-g_m(0)\Big]\Big[\frac{f_{n+1-m}(k)}{f_n(k)}-1\Big].
    \enarray
The term $X_1$ is bounded using Lemma~\ref{lem-pibds}(iv)
with $\epsilon' \in (\delta, \epsilon \wedge 1 \wedge \frac{d-4}{2})$, and
using the
fact that $a(k)\leq \gamma (n+1)^{-1} \log{(n+1)}$, by
\eq
    |X_1|  \leq  K_4' \beta
    a(k)^{1+\epsilon'}\sum_{m=2}^{n+1}
    \frac{1}{m^{(d-2-2\epsilon')/2}} \leq C K_4' \beta
    a(k)^{1+\epsilon'} \leq \frac{CK_4' \beta a(k)}{(n+1)^{\delta}}.
\label{I}
\en

For $X_2$, we first apply Lemma~\ref{lem-pibds}(ii,iv), with
$\epsilon' =0$, to obtain
\eq
    |g_m(k)-g_m(0)| \leq 2 K_4' \beta a(k) m^{-(d-2)/2}.
\en
Applying \refeq{ratio1} then gives
\eq
\lbeq{X2bd}
    |X_2| \leq CK_4' \beta a(k)
    \sum_{m=2}^{n+1} \frac{1}{m^{(d-2)/2}}
    \left(
    (m-1) a(k)(n+1)^{\gamma q}
    + \frac{K_3\beta}{(n+2-m)^{(d-4)/2}}
    \right).
\en
By an elementary estimate, the contribution from the
second term on the right side
is bounded above by $CK_3K_4'\beta^2 a(k) (n+1)^{-(d-4)/2}$.
The first term is bounded above by
\eq
    CK_4' \beta a(k) (n+1)^{0 \vee (6-d)/2} (n+1)^{\gamma q - 1} \log (n+1)
\en
(a harmless factor $\log (n+1)$ should appear for $d=6$).
Since we may choose $q$ to be as close to $1$ as desired, by \refeq{agddef}
this is bounded above by $CK_4' \beta a(k)(n+1)^{-\delta}$.
With (\ref{I}), this proves the bound on $X$ in \refeq{rpk0suf}.

\medskip \noindent {\em Bound on $Y$}.
By \refeq{fs},
    \eq
    \lbeq{AkA0ratio}
    \frac{f_{n+1-m}(k)}{f_n(k)} =
    \frac{f_{n+1-m}(0)}{f_n(0)} \prod_{j=n+2-m}^n
    [1 -v_j a(k)+s_j(k)]^{-1}  .
    \en
Recalling the definition of $R_{m,n}(k)$
in \refeq{Rdef}, we can therefore decompose $Y$ as $Y=Y_1+Y_2+Y_3+Y_4$ with
    \eqarray
    Y_1 & = &  \sum_{m=2}^{n+1} g_m(0)
    \frac{f_{n+1-m}(0)}{f_n(0)}
    R_{m,n}(k),  \\
    Y_2 & = &  \sum_{m=2}^{n+1} g_m(0)
    \frac{f_{n+1-m}(0)}{f_n(0)}
    \sum_{j=n+2-m}^n \left[ (v_j-v_{n+1})a(k)
    - s_j(k) \right], \\
    Y_3 & = &  \sum_{m=2}^{n+1} g_m(0)
    \left[ \frac{f_{n+1-m}(0)}{f_n(0)} - 1 \right]
    (m-1)v_{n+1} a(k), \\
    Y_4 & = &  \sum_{m=2}^{n+1} g_m(0)
    \left[ \frac{f_{n+1-m}(0)}{f_n(0)} - 1 \right].
    \enarray
Then
    \eq
    Y(0) = Y_4  \quad \mbox{ and } \quad
    Y(k)-Y(0) = Y_1+Y_2+Y_3.
    \en

For $Y_1$, we use Lemma~\ref{lem-pibds}, \refeq{ratio0} and \refeq{Rbd}
to obtain
    \eq
    |Y_1| \leq  CK_4'\beta  a(k)^2 (n+1)^{\gamma q}
    \sum\limits_{m=2}^{n} \frac{1}{m^{(d-4)/2}}
    \lbeq{III2z}.
    \en
As in the analysis of the first term of \refeq{X2bd}, we therefore have
    \eq
    \lbeq{III2a}
    |Y_1| \leq  \frac{CK_4' \beta a(k)}{(n+1)^{\delta}}.
    \en
For $Y_2$, we use Lemma~\ref{lem-pibds}, \refeq{ratio0}, (H2), \refeq{sbd}
and Lemma~\ref{lem-conv}
to obtain
    \eq
    \lbeq{II2bd}
    |Y_2| \leq \sum_{m=2}^{n+1} \frac{ K_4' \beta}{m^{d/2}}
    C
    \sum_{j=n+2-m}^n \left[ \frac{K_2 \beta a(k)}{j^{(d-4)/2}}
    + \frac{K_3 \beta a(k)}{j^{\delta}} \right]
    \leq \frac{CK_4' (K_2+K_3) \beta^2 a(k)}{(n+1)^{\delta}}.
    \en
The term $Y_3$ obeys
    \eq
    \lbeq{II3bd}
    |Y_3| \leq
    \sum_{m=2}^{n+1} \frac{K_4'\beta}{m^{(d-2)/2}}
    \frac{CK_3 \beta }{(n+2-m)^{(d-4)/2}}  a(k)
    \leq \frac{CK_4' K_3 \beta^2 a(k)}{(n+1)^{(d-4)/2}},
    \en
where we used Lemma~\ref{lem-pibds}, \refeq{ratio1.a},
\refeq{vnear1}, and an elementary convolution bound.   This proves
the bound on $|Y(k)-Y(0)|$ of \refeq{rpk0suf}, if $\beta$ is
sufficiently small.

We bound $Y_4$ in a similar fashion, using Lemma~\ref{lem-conv}
and the intermediate bound of \refeq{ratio1.a} to obtain
    \eq
    \label{Vest}
    |Y_4|  \leq
    \sum_{m=2}^{n+1} \frac{ K_4'\beta}{m^{d/2}}
    \sum_{j = n+2-m}^n
    \frac{CK_3 \beta}{j^{(d-2)/2}}
    \leq \frac{CK_4'K_3 \beta^2 }{(n+1)^{(d-2)/2}}.
    \en
Taking $\beta$ small then gives the bound on $Y(0)$ of
\refeq{rp0suf}.

\medskip \noindent {\em Bound on $Z$}.
We decompose $Z$ as
\eq
\lbeq{IIIsplit}
    Z = \frac{e_{n+1}(0)}{f_n(0)}
    +\frac{1}{f_n(0)}
    \left[ e_{n+1}(k) - e_{n+1}(0) \right]
    + \frac{e_{n+1}(k)}{f_n(0)}\left[\frac{f_n(0)}{f_n(k)} - 1\right]
    = Z_1+Z_2+Z_3.
\en
Then
        \eq
        Z(0)=Z_1 \quad \mbox{ and } \quad Z(k)-Z(0) = Z_2+Z_3.
        \en
Using Lemma~\ref{lem-pibds}(v,vi), and \refeq{ratio0}
with $m=n+1$, we
obtain
    \eq
    |Z_1| \leq CK_4' \beta (n+1)^{-d/2} \quad \mbox{ and }
    \quad
    |Z_2| \leq CK_4' \beta a(k)(n+1)^{-(d-2)/2}.
    \en
Also, by \refeq{ratio0} and \refeq{ratio2}, we have
\eq
     |Z_3|
     \leq   CK_4' \beta (n+1)^{-d/2}
     a(k) (n+1)^{1 + \gamma q}
     \leq CK_4' \beta a(k)(n+1)^{-(1+\delta)}.
\en

%

This completes the proof of \refeq{rbds}, and hence completes
the advancement of (H3) to $n+1$.

\subsection{Advancement of (H4)}
\label{sec-advH3}\label{sec-advH35}

In this section, we fix $a(k) > \gamma (n+1)^{-1} \log (n+1)$. To
advance (H4) to $j=n+1$, we first recall the definitions of
$b_{n+1}$, $\zeta_{n+1}$ and $X_1$ from (\ref{Delta_n}),
\refeq{zetadef} and \refeq{X1def}.  After some algebra,
\refeq{fkrec} can be rewritten as
    \eq
    f_{n+1}(k) =  f_n(k)\Big( 1~- a(k)b_{n+1}
    + X_1 + \zeta_{n+1}\Big) + W +e_{n+1}(k),
    \label{the eq (H3-H4)}
    \en
with
    \eq
%
    \lbeq{H3IIdef}
    W   = \sum_{m=2}^{n+1} g_m(k) \left[f_{n+1-m}(k)-f_n(k)\right].
    \en

We already have estimates for most of the relevant terms. By
Lemma~\ref{zetan}, $|\zeta_{n+1}| \leq CK_1
\beta(n+1)^{-(d-2)/2}$. By (\ref{I}), $|X_1| \leq CK_4' \beta
a(k)^{1+\epsilon'}$, for any $\epsilon' \in(\delta, \epsilon \wedge
1 \wedge \frac{d-4}{2})$. By
Lemma~\ref{lem-pibds}(v), $|e_{n+1}(k)|\leq K_4' \beta
(n+1)^{-d/2}$. It remains to estimate $W$.  We will show below that $W$
obeys the bound
    \eq
    \lbeq{Wbd}
    |W|
    \leq \frac{CK_4' \beta}{a(k)^{1+\rho}(n+1)^{d/2}}
    (1+K_3\beta +K_5).
    \en
Before proving \refeq{Wbd}, we will first show that it is sufficient
for the advancement of (H4).

In preparation for this, we first note that it suffices to consider only
large $n$.   In fact, since
$|f_n(k;z)|$ is bounded uniformly in $k$ and in $z$ in a
compact set by Assumption~S, and since $a(k)\leq 2$, it is clear
that both inequalities of (H4) hold for all $n\leq N$, if we
choose $K_4$ and $K_5$ large enough (depending on $N$).
We therefore assume in the following
that $n \geq N$ with $N$ large.

Also, care is required to invoke (H3) or (H4), as applicable, in
estimating the factor $f_n(k)$ of (\ref{the eq (H3-H4)}).  Given $k$,
(H3) should be used for the value $n$ for which
$\gamma (n+1)^{-1} \log (n+1) < a(k) \leq \gamma n^{-1} \log n$
((H4) should be used for larger $n$).
We will now show that, as anticipated in the discussion of (H4)
in Section~\ref{sec-mot}, the bound of (H3) actually implies the
bound of (H4) in this case.  To see this, we use Lemma~\ref{lem-cA}
to see that there are $q,q'$ arbitrarily close to $1$ such that
\eq
\lbeq{H3toH4}
    |f_n(k)| \leq Ce^{-qa(k)n} \leq \frac{C}{(n+1)^{q\gamma n/(n+1)}}
    \leq \frac{C}{n^{q' \gamma}}
    \leq \frac{C}{n^{d/2}}\frac{n^{2+\rho}}{n^{q'\gamma +2+\rho -d/2}}
    \leq \frac{C}{n^{d/2}a(k)^{2+\rho}},
\en
where we used the fact that $\gamma + 2 +\rho -d/2 >0$ by \refeq{agddef}.
Thus, taking $K_4 \gg 1$,
we may use the bound of (H4) also for the value of $n$ to which
(H3) nominally applies.  We will do so in what follows, without further comment.

\medskip\noindent {\em Advancement of the second bound of (H4)
assuming \refeq{Wbd}}.
To advance the second estimate in (H4),
we use (\ref{the eq (H3-H4)}), (H4), and the bounds found
above, to obtain
    \eqarray
    \Big|f_{n+1}(k) -  f_n(k)\Big|
    &\leq & \big|f_n(k)\big|~
    \big|-a(k)b_{n+1}+ X_1 + \zeta_{n+1}\big| + |W|+|e_{n+1}(k)| \nonumber \\
    & \leq & \frac{K_4}{n^{d/2}a(k)^{2+\rho}}
    \left(a(k)b_{n+1}+ C K_4' \beta a(k)^{1+\epsilon'}
    + \frac{CK_1\beta }{(n+1)^{(d-2)/2}}
    \right)\nn\\
    \lbeq{H3sec}
    &&\quad+ \frac{CK_4' \beta (1+K_3\beta +K_5)}{(n+1)^{d/2}a(k)^{1+\rho}}
    + \frac{ K_4' \beta}{(n+1)^{d/2}}.
    \enarray
Since $b_{n+1}=1+\Ocal(\beta)$ by \refeq{bnear1}, and
since $(n+1)^{-(d-2)/2} < [a(k) /\gamma\log (n+1)]^{(d-2)/2} \leq C a(k)$, the
second estimate in (H4) follows for $n+1$ provided $K_5 \gg K_4$
and $\beta$ is sufficiently small.

\medskip\noindent {\em Advancement of the first bound of (H4) assuming
\refeq{Wbd}}.
To advance the first estimate of (H4), we argue as in
\refeq{H3sec} to obtain
    \eqarray
    \big|f_{n+1}(k)\big|  &\leq &
    \big|f_n(k)\big|~\Big|1~-a(k)b_{n+1}+ X_1 +\zeta_{n+1}\Big|
    + |W| +|e_{n+1}(k)| \nonumber\\
    &\leq & \frac{K_4}{n^{d/2}a(k)^{2+\rho}}
    \left( |1~-a(k)b_{n+1}|
    + C K_4' \beta a(k)^{1+\epsilon'}
    + \frac{CK_1\beta }{(n+1)^{(d-2)/2}}
    \right)\nn\\
    &&\quad+ \frac{CK_4' \beta (1+K_3\beta +K_5)}{(n+1)^{d/2}a(k)^{1+\rho}}
    + \frac{ K_4' \beta}{(n+1)^{d/2}}.
    \lbeq{H3bd1}
    \enarray
We need to argue that the right-hand side is no larger than $K_4
(n+1)^{-d/2} a(k)^{-2-\rho}$. To achieve this, we will use
separate arguments for $a(k) \leq \frac 12$ and $a(k) > \frac 12$.
These arguments will be valid only when $n$ is large enough.

Suppose that $a(k) \leq \frac 12$. Since $b_{n+1} = 1+
\Ocal(\beta)$ by \refeq{bnear1}, for $\beta$ sufficiently small we
have
    \eq
    1~-b_{n+1}a(k) \geq 0.
    \en
Hence, the absolute value signs on the right side of \refeq{H3bd1}
may be removed.
Therefore, to obtain the first estimate of (H4) for $n+1$,
it now suffices to show that
    \eq
    \lbeq{H3bd2}
    1~-ca(k)
    +\frac{CK_1\beta}{(n+1)^{(d-2)/2}}
    \leq \frac{n^{d/2}}{(n+1)^{d/2}},
    \en
for $c$ within order $\beta$ of 1.  The term $ca(k)$ has been
introduced to absorb $b_{n+1}a(k)$, the order $\beta$ term in
\refeq{H3bd1} involving $a(k)^{1+\epsilon'}$, and the last two
terms of \refeq{H3bd1}.  However,
$a(k) > \gamma (n+1)^{-1} \log(n+1)$.  From this, it can be
seen that \refeq{H3bd2}
holds for $n$ sufficiently large and $\beta$ sufficiently small.

Suppose, on the other hand, that $a(k) > \frac 12$. By
\refeq{Dbound3}, there is a positive $\eta$, which we
may assume lies in $(0,\frac{1}{2})$, such that $-1+\eta <1-a(k)<
\frac 12$. Therefore $|1-a(k)| \leq 1-\eta$ and
    \eq
    |1~-b_{n+1}a(k)|
    \leq |1-a(k)| + |b_{n+1}-1| \, |a(k)|
    \leq 1-\eta
    + 2|b_{n+1}-1|.
    \en
Hence
    \eq
    |1-a(k)b_{n+1}| + C K_4' \beta a(k)^{1+\epsilon'}
    +\frac{CK_1\beta }{(n+1)^{(d-2)/2}}
    \leq 1-\eta + C(K_1+K_4')\beta,
    \en
and the right side of \refeq{H3bd1} is at most
    \eqarray
    &&\frac{K_4}{n^{d/2}a(k)^{2+\rho}}
    \left[1-\eta + C(K_1 +K_4')\beta \right] +
    \frac{CK_4'(1+K_3\beta + K_5) \beta}{(n+1)^{d/2}a(k)^{2+\rho}}\nonumber\\
    &&\quad\leq  \frac{K_4}{n^{d/2}a(k)^{2+\rho}}
    \left[1-\eta + C(K_5K_4'+K_1)\beta \right].
    \enarray
This is less than $K_4 (n+1)^{-d/2}a(k)^{-2-\rho}$ if $n$ is
large and $\beta$ is sufficiently small.

This advances the
first bound in (H4), assuming \refeq{Wbd}.

\medskip\noindent {\em Bound on $W$}.
We now obtain the bound \refeq{Wbd} on $W$.  As a first step,
we rewrite $W$ as
    \eq
    \lbeq{Wdef}
    W  =\sum_{j=0}^{n-1} g_{n+1-j}(k)\sum_{l=j+1}^n
    [f_{l-1}(k)-f_l(k)].
    \en
Let
    \eq
    m(k) = \max \{ l\in \{2,\ldots, n\} : a(k) \leq \gamma l^{-1} \log l \}.
    \en
For $l \leq m(k)$, $f_l$ is in the domain of (H3), while for $l>m(k)$,
$f_l$ is
in the domain of (H4).  By hypothesis, $a(k) > \gamma (n+1)^{-1} \log(n+1)$.
We divide the sum over $l$ into two parts,
corresponding respectively to $l \leq m(k)$ and $l> m(k)$,
yielding $W=W_1 + W_2$,
By Lemma~\ref{lem-pibds}(i),
    \eqarray
    |W_1|  & \leq &\sum_{j=0}^{m(k)} \frac{ K_4' \beta}{(n+1-j)^{d/2}}
    \sum_{l=j+1}^{m(k)} |f_{l-1}(k)-f_l(k)|
    \\
    |W_2|  & \leq &\sum_{j=0}^{n-1}\frac{ K_4' \beta}{(n+1-j)^{d/2}}
    \sum_{l=(m(k) \vee j)+1}^n |f_{l-1}(k)-f_l(k)|.
    \enarray
The term $W_2$ is easy, since by (H4) and Lemma \ref{lem-conv} we have
    \eq
    \label{sum1}
    |W_2| \leq \sum\limits_{j=0}^{n-1} \frac{
    K_4' \beta}{(n+1-j)^{d/2}} \sum\limits_{l=j+1}^n
    \frac{K_5}{a(k)^{1+\rho} \; l^{d/2}} \leq \frac{C K_5 K_4'
    \beta}{ a(k)^{1+\rho} (n+1)^{d/2}}.
    \en



For $W_{1}$, we have the estimate
\eq
\lbeq{W1pbd}
    |W_1| \leq \sum_{j=0}^{m(k)} \frac{ K_4' \beta}{(n+1-j)^{d/2}}
    \sum_{l=j+1}^{m(k)} |f_{l-1}(k)-f_l(k)|.
\en
For $1 \leq l \leq m(k)$, it follows from Lemma~\ref{lem-cA}
and \refeq{ratio3} that
    \eq
    \lbeq{fldiff}
    |f_{l-1}(k)-f_l(k)|
    \leq
    C e^{- q a(k) l}\left( a(k) + \frac{K_3 \beta}{l^{(d-2)/2}}\right),
    \en
with $q =1-\Ocal(\beta)$.
We fix a small $r > 0$, and bound the summation over $j$
in \refeq{W1pbd} by summing separately over $j$ in the ranges $1
\leq j \leq (1-r)n$ and $(1-r)n \leq j \leq m(k)$ (the
latter range may be empty).  We denote the contributions from
these two sums by $W_{1,1}$ and $W_{1,2}$ respectively.

To estimate $W_{1,1}$, we will make use of the bound
\eq
    \sum_{l=j+1}^\infty e^{- q a(k) l}  l^{-p}
    \leq
    \left\{
    \begin{array}{ll}
    C e^{- q a(k) j} a(k)^{p-1} & (0 \leq p < 1) \\
    C e^{- q a(k) j} & (p>1),
    \end{array}
    \right.
\en
which can be proved by approximating the sum by an integral.
With \refeq{W1pbd} and \refeq{fldiff}, this gives
    \eqarray
    |W_{1,1}|
    & \leq & \frac{CK_4' \beta}{(n+1)^{d/2}}
    \sum_{j=0}^{(1-r)n}
    e^{- q a(k) j}
    \left( 1 + K_3 \beta\right)
    \nn \\
    \lbeq{W11p}
    & \leq & \frac{CK_4' \beta}{(n+1)^{d/2}} \frac{1+K_3\beta}{a(k)}
    \leq \frac{CK_4' \beta}{(n+1)^{d/2}} \frac{1+K_3\beta}{a(k)^{1+\rho}}.
    \enarray

For $W_{1,2}$, we have
    \eq
    |W_{1,2}| \leq   \sum_{j=(1-r)n}^{m(k)}
    \frac{CK_4'\beta }{(n+1-j)^{d/2}}
    \sum_{l=j+1}^{m(k)}
     e^{-q a(k) l}\left(a(k) + \frac{K_3 \beta}{l^{(d-2)/2}}\right).
    \en
Since $l$ and $m(k)$ are comparable and large, it follows as in
\refeq{H3toH4} that 
\eq
    e^{-q a(k) l}
    \left(a(k) + \frac{K_3 \beta}{l^{(d-2)/2}}\right)
    \leq \frac{C}{a(k)^{2+\rho} l^{d/2}}
    \left(a(k) + \frac{K_3 \beta}{l^{(d-2)/2}}\right)
    \leq \frac{C(1+K_3\beta)}{a(k)^{1+\rho}l^{d/2}} .
\en
Hence, by Lemma~\ref{lem-conv},
    \eq
    \lbeq{W12p}
    |W_{1,2}| \leq \frac{C(1+K_3\beta)K_4' \beta}{a(k)^{1+\rho}}
    \sum_{j=(1-r)n}^{m(k)} \frac{1}{(n+1-j)^{d/2}}
    \sum_{l=j+1}^{m(k)} \frac{1}{l^{d/2}}
    \leq \frac{C(1+K_3\beta)K_4' \beta}{a(k)^{1+\rho}(n+1)^{d/2}}.
    \en

Summarising, by \refeq{W11p}, \refeq{W12p},
and (\ref{sum1}), we have
    \eq
    |W| \leq |W_{1,1}| + |W_{1,2}|  +|W_2|
    \leq \frac{CK_4' \beta}{a(k)^{1+\rho}(n+1)^{d/2}}
    (1+K_3\beta +K_5),
    \en
which proves \refeq{Wbd}.

\section{Proof of the main results}
\label{sec-pf} \setcounter{equation}{0}

As a consequence of the completed induction, it follows from
Lemma~\ref{lem-In} that $I_1 \supset I_2 \supset I_3 \supset
\cdots$, so $\cap_{n=1}^\infty I_n$ consists of a single point
$z=z_c$.  Since $z_0=1$, it follows from (H1) that $z_c=1+\Ocal(\beta)$.
We fix $z=z_c$ throughout this section.  The constant $A$ is
defined by $A = \prod_{i=1}^\infty [1+
r_i(0)] = 1+ {\cal O}(\beta)$.   To define $v$, we note that
it follows from (H2) that there is a
constant $v = 1+ {\cal O}(\beta)$ such that
        \eq
        \lbeq{Delta_nlimit}
        |v_n(z_c) - v| \leq {\cal O}( \beta n^{-(d-4)/2}).
        \en

\subsection{Proof of Theorem \ref{thm-1}}

\smallskip \noindent {\em Proof of Theorem \ref{thm-1}(a).}
By (H3),
        \eq
        \lbeq{hat{tau}_nlimit}
        |f_n(0;z_c) - A|
        = \prod_{i=1}^n [1+r_i(0)]
        \big| 1 - \prod_{i=n+1}^\infty [1+r_i(0)]\big|
        \leq {\cal O}(\beta n^{-(d-4)/2}).
        \en
Suppose $k$ is such that $a(k/\sqrt{Dn}) \leq \gamma n^{-1} \log
n$, so that (H3) applies.  Here, we use the $\gamma$ of \refeq{agddef}.
By \refeq{momentD},
$a(k) = \sigma^2k^2/2d  + {\cal O}(k^{2+2\epsilon})$ with $\epsilon > \delta$,
where we now allow constants in error terms to depend on $L$.
Using this, together with \refeq{fs}--\refeq{sbd},
\refeq{Delta_nlimit}, and  $\delta <
1 \wedge \frac{d-4}{2}$, we obtain
        \eqarray
        \frac{f_n(k/\sqrt{v\sigma^2 n};z_c)}{f_n(0;z_c)}
        & = &
        \prod_{i=1}^n \left[1- v_i a\big(\frac{k}{\sqrt{v\sigma^2 n}}\big)
        + {\cal O}(\beta a\big(\frac{k}{\sqrt{v\sigma^2 n}}\big) i^{-\delta})
        \right]
        \nnb
        & = &
        \lbeq{11cpf}
        e^{-k^2/2d}[1+{\cal O}(k^{2+2\epsilon} n^{-\epsilon} )
        + {\cal O}(k^2 n^{- \delta})].
        \enarray
With \refeq{hat{tau}_nlimit}, this gives the desired result.

\medskip \noindent {\em Proof of Theorem \ref{thm-1}(b).}
Since $\delta<1 \wedge
\frac{d-4}{4}$, it follows from \refeq{1.2b2}--\refeq{1.2b3}
and \refeq{Delta_nlimit}--\refeq{hat{tau}_nlimit} that
        \eq
        \frac{\nabla^2f_n(0;z_c)}{f_n(0;z_c)}
        = -v\sigma^2  n [1+{\cal O}(\beta n^{-\delta})].
        \en

\medskip \noindent {\em Proof of Theorem \ref{thm-1}(c).} The claim is
immediate from Lemma~\ref{lem-L1norm}, which is now known to hold for all $n$.

\subsection{Proof of Theorem~\protect\ref{thm-zcpc}}

By Theorem~\ref{thm-1}(a), $\chi(z_c)=\infty$.  Therefore $z_c \geq p_c$.
We need to rule out the possibility that $z_c>p_c$.
Theorem~\ref{thm-1} also
gives \refeq{fbds} at $z=z_c$.  By assumption, the series
\eq
    G(z) = \sum_{m=2}^\infty g_m(0;z),
    \quad
    E(z) = \sum_{m=2}^\infty e_m(0;z)
\en
therefore both converge and are ${\cal O}(\beta)$ uniformly in $z \leq z_c$.
For $z<p_c$, the basic
recursion relation \refeq{fkrec} gives
\eq
    \chi(z) = 1 + z\chi(z) + G(z)\chi(z) +E(z),
\en
and hence
\eq
\lbeq{chiEG}
    \chi(z) = \frac{1+E(z)}{1-z-G(z)}, \quad (z<p_c).
\en
Since $E(z) = {\cal O}(\beta)$ and
$\lim_{z \uparrow p_c}\chi(z) = \infty$, it follows from \refeq{chiEG} that
\eq
\lbeq{pc0}
    1 - p_c - G(p_c) = 0.
\en
By (\ref{z}), \refeq{pc0} holds also when $p_c$ is replaced by $z_c$.
If $p_c \neq z_c$, then it follows from the mean-value theorem that
\eq
    z_c - p_c = G(p_c) - G(z_c) = -(z_c-p_c) \sum_{m=2}^\infty
    \partial_z g_m(0;t)
\en
for some $t \in (p_c,z_c)$.  However, by a bound of Assumption~G, the sum
on the right side is ${\cal O}(\beta)$ uniformly in $t \leq z_c$.
This is a contradiction, so we
conclude that $z_c=p_c$.

\subsection{Proof of Theorem~\protect\ref{thm-3}}

We begin by noting that
\eq
    \frac{1}{(2R+1)^d}
    \sum_{y \in C_R( \lfloor x \sqrt{v\sigma^2 n} \rfloor )} p_n(y)
    =
    (Q_R * p_n)(\lfloor x\sqrt{n}\rfloor ),
\en
where $Q_R(x)  = (2R+1)^{-d}$ for $x\in C_R(0 )$, and otherwise equals zero.
Therefore,
    \eq
    \lbeq{int1}
    \frac{1}{(2R+1)^d} \sum_{y\in C_R(\lfloor x\sqrt{v\sigma^2 n}\rfloor )}
    p_n(y) =
    \int_{[-\pi,\pi]^d}
    \frac{d^d k}{(2\pi)^d} e^{-ik\cdot \lfloor x\sqrt{v\sigma^2 n}\rfloor}
    \hat{Q}_R(k)f_n(k;z_c).
    \en
Before proceeding with the proof, we first derive some relevant properties of
$Q_R$.

By definition, $|\hat{Q}_R(k)| \leq 1$.
The Fourier transform $\hat{Q}_R(k)$ is given in terms of the
Dirichlet kernel
\eq
    M(t) = \sum_{j=-R}^R e^{ijt} =
    \frac{\sin [(2R+1)t/2]}{\sin (t/2)}
\en
by
    \eq\lbeq{DReq}
    \hat{Q}_R(k)=\prod_{i=1}^d \frac{M(k_i)}{2R+1}.
    \en
Since $|M(t)|$ is bounded above by both $2R+1$ and $|\sin (t/2)|^{-1}$, we
have
    \eq \lbeq{Ombd}
    | \hat{Q}_R(k) |
    \leq
    \frac{1}{(2R+1)\sin (\|k\|_\infty /2)}
    \leq
    \frac{\pi}{2 ||k||_{\infty}R}.
    \en
Let $\kappa \in (0,1)$.  When the right side of \refeq{Ombd} is less than $1$,
the bound is degraded by raising the right side to the power $\kappa$.
On the other hand, when the right side is greater than $1$, we may
still raise it to the power $\kappa$ and have a correct bound, since the left
side is at most $1$.  Therefore, for every $\kappa \in (0,1)$, we have
    \eq
    \label{DRbound}
    |\hat{Q}_R(k)| \leq \Big[\frac{\pi}{2 ||k||_{\infty}R}\Big]^\kappa.
    \en

We divide the domain of integration
$[-\pi,\pi]^d$ of \refeq{int1} into a small-$k$ region
$S_n = \{k \in [-\pi,\pi]^d : a(k)\leq \gamma n^{-1} \log
n\}$ and a large-$k$ region
$L_n = \{k \in [-\pi,\pi]^d : a(k)> \gamma n^{-1} \log
n\}$.  By Theorem~\ref{thm-1}(a), the integral over the small-$k$ region
equals
    \eq
    A\int_{S_n} \frac{d^d k}{(2\pi)^d}
    e^{-ik\cdot \lfloor x\sqrt{v\sigma^2 n}\rfloor}
    e^{-nDk^2/2d}[1+{\cal O}(k^2 n^{1-\delta})+{\cal O}(k^2 R^2)] ,
    \en
where we allow constants in error terms to depend on $L$.
For the error terms, we bound the complex exponential by $1$ and obtain a
factor $n^{-d/2}$ from the integration.
For the leading term, we extend the integration domain
to $\Rd$ and perform the integral exactly.  The error incurred by extending
the integration domain is at most $O(n^{-d/2-\gamma^-})$, where $\gamma^-$
is any positive number less than $\gamma$.  Therefore, the small-$k$ region
gives
\eq
    A\left(\frac{d}{2\pi v\sigma^2 n }\right)^{d/2} \left[e^{-dx^2/2}
    + \Ocal(n^{-\gamma^-})+ {\cal O}(n^{-\delta})
    +{\cal O}(R^2 n^{-1}) \right].
\en

The integral over the large-$k$ region can be bounded using (H4),
(\ref{DRbound}) and \refeq{Dbound1} by
    \eq
    K_4 n^{-d/2} \int_{L_n} \frac{d^d k}{(2\pi)^d}
    |\hat{Q}_R(k)| a(k)^{-2-\rho} \leq \Ocal(n^{-d/2}R^{-\kappa})
    \int_{L_n} d^d k \,
    k^{-4-2\rho - \kappa}.
    \en
The integral is bounded provided
$\kappa < d-4 + 2\rho$.  By \refeq{agddef}, we may choose a positive $\kappa$
obeying this bound.

Therefore, the left side of \refeq{int1} is given by
\eq
\lbeq{LCLTerr}
    A\left(\frac{d}{2\pi v\sigma^2 n}\right)^{d/2} \left[e^{-dx^2/2}
    + \Ocal(R^{-\kappa})
    + \Ocal(n^{-\gamma^-})+ {\cal O}(n^{-\delta})
    +{\cal O}(R^2 n^{-1}) \right].
\en
For $x^2$ less than a sufficiently small multiple of $\log R$,
this has the desired asymptotic form, provided $R=R_n \to \infty$
as $n \to \infty$, with $R_n = {o}(n^{1/2})$.
\qed

\subsection{Formulae for $z_c$, $A$, $v$}
\label{sec-zAD}

The identity (\ref{z})
follows after we let $n \to \infty$ in \refeq{zetadef}, using
Lemma~\ref{zetan}.

To determine $A$, we use a summation argument. Let $\chi_n =
\sum_{k=0}^n f_k(0)$.  By \refeq{fkrec},
    \eqarray
    \chi_n &=& 1 + \sum_{j=1}^n f_j(0)
    = 1 + \sum_{j=1}^n \sum_{m=1}^{j}g_m(0) f_{j-m}(0) +  \sum_{j=1}^n e_j(0)
    \nonumber\\
    &=& 1 + z\chi_{n-1}
    +\sum_{m=2}^{n}g_m(0)\chi_{n-m}+ \sum_{m=1}^ne_m(0).
    \label{S_nsum}
    \enarray
Using \refeq{zetadef} to rewrite $z$, this gives
    \eq
    f_n(0) = \chi_n - \chi_{n-1}
    = 1+ \zeta_n \chi_{n-1} - \sum_{m=2}^{n}
    g_m(0)(\chi_{n-1}-\chi_{n-m}) + \sum_{m=1}^ne_m(0).
    \en
By Theorem~\ref{thm-1}(a), $\chi_n \sim nA$ as $n \to \infty$.
Therefore, using Lemma~\ref{zetan} to bound the $\zeta_n$ term,
taking the limit $n \to \infty$ in the above equation gives
    \eq
    A = 1 - A \sum_{m=2}^\infty (m-1)g_m(0) + \sum_{m=1}^\infty e_m(0).
    \en
With (\ref{z}), this gives (\ref{A}).

To prove (\ref{D}), we use \refeq{Delta_nlimit}, (\ref{Delta_n}) and
Lemma~\ref{lem-pibds} to obtain
    \eq
    v = \lim_{n \to \infty}  v_n
    = \frac{-\sigma^{-2}\sum_{m=2}^\infty \nabla^2 g_m(0)}
    {1 + \sum_{m=2}^\infty (m-1) g_m(0)}.
    \en
Equation~(\ref{D}) then follows, once we rewrite the denominator
using (\ref{z}).

\renewcommand{\thesection}{\Alph{section}}
\setcounter{section}{0}

\section{Appendix: Example of $D_L$}
\label{sec-app}
\setcounter{equation}{0}

In this Appendix, we verify that the function $D_L$ defined in \refeq{Dbasic2}
obeys the requirements of Assumption~D.
The bounds of  \refeq{momentD} and \refeq{beta, sigmadef}
follow easily, and clearly $\hat{D}_L(0)=1$.
We therefore concentrate on
establishing \refeq{Dbound1} and \refeq{Dbound2}--\refeq{Dbound3}.
Before proceeding, we first note that it is sufficient to prove the lower
bounds of \refeq{Dbound1} and \refeq{Dbound2} respectively for
$\|k\|_\infty \leq bL^{-1}$ and $\|k\|_\infty \geq bL^{-1}$,
for any small positive $b$, since together these bounds imply the
corresponding statements with $b=1$.  We will prove this modified form
of the bounds.

For the upper bound of \refeq{Dbound1}, we use only the symmetry of
$D_L$ and the inequality $1-\cos{t}\leq \frac 12 t^2$ to obtain
\eq
    0\leq a_L(k) = \sum_{x\in \Z^d} [1-\cos{(k\cdot x)}] D_L(x)
    \leq \frac12 \sum_{x\in \Z^d}(k\cdot x)^2 D_L(x)=\frac12
    k^2 \sigma^2
    \quad
    (k \in [-\pi,\pi]^d).
\en
For a lower bound, we use
$1-\cos{t}\geq \frac 12 t^2 - \mbox{const.} t^{2+2\epsilon}$,
\refeq{momentD} and \refeq{Dbasic2} to conclude that there is a small
positive $b$ such that
\eq
    a_L(k) \geq
    \frac12 k^2 \sigma^2 - \mbox{const.} k^{2+2\epsilon} L^{2+2\epsilon}
    \geq \mbox{const.} L^2 k^2
    \quad
    (\|k\|_\infty \leq bL^{-1}).
\en

For \refeq{Dbound2}--\refeq{Dbound3}, it suffices to
show that $\hat{D}_L(k)$ is bounded away from $1$ for $\|k\|_\infty \geq bL^{-1}$
and bounded away from $-1$ for $k \in [-\pi,\pi]^d$.
Let $\tilde{h}(k) = \int_{\Rd} h(x) e^{ik\cdot x} d^dx$ denote
the Fourier integral transform of $h$.
Then it suffices to show that
\eq
\lbeq{uc}
    |\hat{D}_L(k/L) - \tilde{h}(k)| \to 0
\en
uniformly in $k \in [-L\pi,L\pi]^d$, since the characteristic function
$\tilde{h}$ of a
piecewise continuous probability density is bounded away from $-1$, and is
bounded away from $1$ on the complement of any compact neighbourhood of $0$.
In the remainder of the proof, we prove \refeq{uc}.

Let $X_L$ be a discrete random variable with probability mass function
\eq
    \Pbold (X_L = z ) = D_L(Lz)
    = \frac{h(z)}{\sum_{w\in L^{-1}\Zd} h(w) }
    \quad (z \in L^{-1}\Zd).
\en
Then $X_L$ has characteristic function $\hat{D}_L(k/L)$.
By definition,
\eq
    \hat{D}_L(k/L)
    = \frac{\sum_{z\in L^{-1}\Zd} h(z) e^{ik \cdot z}}
    {\sum_{z\in L^{-1}\Zd} h(z) }.
\en
The right side involves Riemann sums, and hence
\eq
\lbeq{Dhlim}
    \lim_{L \to \infty} \hat{D}_L(k/L) = \tilde{h}(k)
     \quad \mbox{for $k = {\cal O}(1)$.}
\en
We want to extend \refeq{Dhlim} to uniform convergence for
$k \in [-L\pi,L\pi]^d$.

Given $y \in \Rd$, let $z_y$ denote the closest point in $L^{-1}\Zd$ to $y$,
with a fixed arbitrary rule used to choose from among the closest points
when there is not a unique closest point.
Let $Y_L$ be a continuous random variable with density $h_L$ given by
    \eq
    h_L(y)= \frac{L^dh(z_{y})}{\sum_{z\in L^{-1}\Zd} h(z) }
    \quad (y \in \Rd).
    \en
Thus $Y_L$ has a density which is a piecewise constant analogue of the
probability mass function of $X_L$.  It follows that $Y_L$ has the
same distribution as $X_L+U_L$, where $U_L$ is uniform on
$[-\frac{1}{2L},\frac{1}{2L}]^d$ and independent of $X_L$.
Let $\tilde{h}_L$ and $\tilde{u}_L$ respectively denote the characteristic
functions of $Y_L$ and $U_L$.
Then
    \eq
    \lbeq{Drewrite2}
    \hat{D}(k/L)  = \frac{\tilde{h}_{L}(k)}{\tilde{u}_L(k)}
    \quad (k \in \Rd) .
    \en

We claim that
    \eq
    \lbeq{clim}
    |\tilde{h}_{L}(k)-\tilde{h}(k)| \rightarrow 0, \quad\mbox{
    uniformly in $k \in \Rd$.}
    \en
To see this, we fix $\epsilon >0$ and choose $R$, independent of $L$,
such that $\int_{\Rd \backslash [-R,R]^d} |h(y)|d^dy < \epsilon$
and $\int_{\Rd \backslash [-R,R]^d} |h_L(y)|d^dy < \epsilon$.
This is possible since
$h$ is integrable, using the monotonicity assumption on $h$.
Since $h$ is almost everywhere continuous, for almost all $y\in \R^d$ we have
    \eq
    \lbeq{hLh}
    \lim_{L\rightarrow \infty} h_L(y) = h(y).
    \en
Therefore, by dominated convergence, $\int_{[-R,R]^d} |h(y) - h_L(y)|d^dy <
\epsilon$ for $L$ sufficiently large.  The claim
\refeq{clim} then follows from
the inequality
$|\tilde{h}_{L}(k)-\tilde{h}(k)| \leq \|h_L - h\|_1$.


By definition,
        \eq
        \lbeq{uLk}
        \tilde{u}_L(k)
        = L^d \int_{[-\frac{1}{2L},\frac{1}{2L}]^d} e^{ik\cdot x} d^dx
        = \prod_{i=1}^d
        \frac{2L\sin{(\frac{k_i}{2L})}}{k_i}.
        \en
Therefore, for any fixed $\alpha >0$, $\tilde{u}_L(k)\rightarrow 1$ uniformly
in $k$ such that $\|k\|_{\infty} \leq L^{1-\alpha}$.  With \refeq{Drewrite2}
and the uniform convergence of \refeq{clim}, this proves that the
convergence in \refeq{uc} holds uniformly in $k$ such that
$\|k\|_{\infty} \leq L^{1-\alpha}$.

It remains to consider the case $\|k\|_{\infty} \geq L^{1-\alpha}$.
By \refeq{uLk}, $\tilde{u}_L(k) \geq (2\pi^{-1})^d$ uniformly in
$k\in [-L\pi, L\pi]^d$.  Also, by the Riemann--Lebesgue lemma,
$\lim_{L \to \infty}
\sup_{\{k: \|k\|_{\infty} \geq L^{1-\alpha}\}}\tilde{h}(k) = 0$.
Therefore, by \refeq{Drewrite2} and \refeq{clim},
\eq
    |\hat{D}(k/L) - \tilde{h}(k)|
    \leq \big(\frac{\pi}{2}\big)^d | \tilde{h}_L(k)| + |\tilde{h}(k)|
    \to 0
\en
uniformly in $k$ such that $\|k\|_{\infty} \geq L^{1-\alpha}$.
This completes the proof of \refeq{Dbound2}--\refeq{Dbound3}.


\section*{Acknowledgements}
This work was supported in part by NSERC of Canada.  The work of both
authors was carried out in part at Microsoft Research and the
Fields Institute.  The work of G.S.\ was also carried out in part
at McMaster University.

\end{document}